\documentclass{article}

\usepackage[ansinew]{inputenc}
\usepackage{theorem}
\usepackage[english]{babel}
\usepackage{graphicx}
\usepackage{amsfonts}
\usepackage{amsmath}

\numberwithin{equation}{section}

\newtheorem {defn}{Definition}[section]
\newtheorem {thm}{Theorem}[section]
\newtheorem {prop}[thm]{Proposition}

\newtheorem {lemma}[thm]{Lemma}
\newtheorem {assump}[defn]{Assumption}
\newtheorem {cor}[thm]{Corollary}

\def\squarebox#1{\hbox to #1{\hfill\vbox to #1{\vfill}}}
\newcommand{\qed}{\hspace*{\fill}
\vbox{\hrule\hbox{\vrule\squarebox{.667em}\vrule}\hrule}\smallskip}

\newcommand{\bbmatrix}[1]{\left[ \begin{array}{cccccccccc} #1 \end{array} \right]}

\title{Overdetermined $2D$ Systems Invariant in One Direction and Their Transfer Functions}
\author{Andrey Melnikov, Victor Vinnikov\\
{\small Ben-Gurion University of the Negev, Israel}
}

\begin{document}
\maketitle
\begin{abstract}
In this work we develop a theory of \textit{Vessels}. This object arises \cite{bib:Vortices, bib:VinnMSRI,bib:Over2Dsys}
in the study of overdetermined $2D$ systems invariant in one of the variables, which are usually called \textit{time invariant}. To each overdetermined time invariant $2D$ systems there is associated a vessel, which 
is a collection of system operators satisfying certain relations and vise versa.
Such an invariance forces the theory of vessels to resemble a constant (classical) $1D$ case \cite{bib:HarmAnal,bib:JordanBrodskii,bib:BallCohen,bib:Inter, bib:BGKD} and 
as a result many notions are naturally redefined and most theorems are reproved in this setting. The notion of transfer function and its connection to
the overdetermined $2D$ time invariant system (and the corresponding vessel) is one of the topics of this work. It is well known 
\cite{bib:BallCohen,bib:BGKD}
that multiplicative structure of a transfer function of a $1D$ system is closely connected to the decomposition of the state space into invariant subspaces of the state operator and we generalize this result to a wider class of functions. This class (denoted by $\boldsymbol {\mathcal I}$)
arises as a class of transfer functions, which intertwine solutions of ODEs with spectral parameters \cite{bib:CoddLev}.
At the end we present solution of factorization problems for finite dimensional case.
\end{abstract}

\tableofcontents

\pagebreak
\section{Introduction}
The theory of two-dimensional ($2D$) overdetermined time-invariant systems has been extensively developed over the last 20 years;
it is closely connected to the theory of commuting operators \cite{bib:TheoryNonComm}, \cite{bib:Over2Dsys}, \cite{bib:VinnMSRI}. 
An overdetermined $2D$ continuous time-invariant linear i/s/o system is of the form
\[
    \Sigma: \left\{ \begin{array}{lll}
    \frac{\partial}{\partial t_1}x(t_1,t_2) = A_1 x(t_1,t_2) + B_1 u(t_1,t_2) \\
    \frac{\partial}{\partial t_2}x(t_1,t_2) = A_2 x(t_1,t_2) + B_2 u(t_1,t_2) \\
    y(t_1,t_2) = D u(t_1,t_2) + C x(t_1,t_2)
    \end{array} \right.
\]
where $(u(t_1,t_2),x(t_1,t_2),y(t_1,t_2))$ is the (input,state,output) triple and all the other symbols denote bounded operator 
on suitable Hilbert spaces. Assuming continuously differentiable inputs, one obtains \cite{bib:Over2Dsys} 
that the state space must be twice differentiable and enjoy the equality of mixed variables, from where algebraic relations
are imposed on the operators of the system, and compatibility conditions (hence {\em overdetermined} system) on the input
$u(t_1,t_2)$ and the output $y(t_1,t_2)$. More precisely, $u(t_1,t_2)$ and $y(t_1,t_2)$ satisfy algebraic equations (live on a curve) 
and an example of relation on the operators is commutativity of $A_1,A_2$.

Using frequency domain analysis a notion of transfer function $S(\lambda,t_2)$ arises
\[ S(\lambda) = D + C(\lambda I - A_1)^{-1} B_1
\]
and the main question is how properties of such transfer functions and properties of the system operators (more precisely invariant subspaces
of $A_1, A_2$) are connected. Many standard structural properties, e.g, controllability, observability, minimality, pairing and adjoint system,
cascade connection, equivalence and standard problems, e.g., pole placement, linear-quadratic-regulator problem $H^\infty$ control
for $1D$ linear systems carry over for this setting.

The study of time varying $1D$ systems has produced a rich theory and the core of all these theories is the Sz.-Nagy–Foias theory of
contractions on a Hilbert space \cite{bib:HarmAnal}. 
The analogue of transfer function in this setting is a lower block-triangular bounded operator. In \cite{bib:AlBallPerez}, for example
there is presented the development of a unified approach to time-varying dissipative linear systems,
non-stationary Lax-Phillips scattering theory, and operator model theory for the infinite family
of contractions. The abstract interpolation problem for the time-varying case is presented both in the de
Branges-Rovnyak model formulation \cite{bib:DeBRModelSpace}, and in its coordinate-free, scattering-theoretic form, and its
application to the time-varying version of the matrix right tangential Nevanlinna-Pick problem is
studied.

Probably, the simplest generalization of time varying $1D$ systems is the study of $2D$ systems invariant in one direction.
There are some works in this direction \cite{bib:Vortices,bib:Gauchman} in different settings. Our main inspiration comes
from the article of M. Liv\v sic \cite{bib:Vortices} and we actually continue this work.
 
So, we suppose that our $2D$ systems are invariant in one of the variables ($t_1$). We will use integrated form for of a system
\cite{bib:Latushkin} and as a result there arises a continuous family of Hilbert spaces and semi-group acting between them.
The invariance of the $2D$ system in one of the variables allows us to perform a partial separation of variables and to define a transfer
function, depending on the corresponding spectral parameter (say $\lambda$), which will, additionally, depend on the second
variable ($t_2$):
\[ S(\lambda,t_2) = D(t_2) + C(t_2)(\lambda I - A_1(t_2))^{-1} B_1(t_2)
\]
A fundamental feature in the study of transfer functions and their factorizations is that algebraic equations are now 
replaced by ODEs with a spectral parameter $\lambda$.

The theory of these systems is interesting by itself, especially since it allows us to use frequency domain
analysis in a time varying framework. It also has important connections with completely integrable nonlinear
PDEs \cite{bib:LaxEquation}: the so-called Lax equation
\[ \frac{d}{dt_2} A_1(t_2) = A_1(t_2) A_2(t_2) - A_2(t_2) A_1(t_2)
\]
appears naturally, and the passage from the input to the output ODE with a spectral parameter is analogous to the B\" acklund transformation.

Let us introduce the detailed description of the contents of this work, highlighting the main results. 
Important notice is that we start the study of $t_1$-invariant $2D$ systems from the more general integrated form, i.e., differential equations are
presented as integral equations and as a result we obtain much weaker assumptions on the operators. Inspiration for doing it in this way comes from
\cite{bib:Latushkin}.
As always the case in this passage, we consider
evolution semi groups acting between continuous set of Hilbert spaces. This semi-group, after appropriate transformation of the Hilbert spaces
and differentiation gives rise to operator $A_2(t_2)$ in section \ref{sec:Vessel}.

There are seven sections in this work. After introduction, at the second section we introduce ovedetermined $2D$ systems and show how basic notions of system theory carry over. Following the ideas in \cite{bib:VinnMSRI} we introduce two notions (local and global) of approximate controllability
and similarly two notions for observability and study their relations. On this basis, we study in section \ref{sec:Vessel} guage quasi-similarity
of (minimal) systems (vessels). One of the interesting results of our work is that there exists a notion of \textit{differential vessel} and
we show that it is always possible to pass from integral to differential form and vice versa. Further, we present the notion of transfer function
and its main properties. As a result we define a class $\boldsymbol{\mathcal I}$ of intertwining functions, which is extensively studied at the last
section \ref{sec:AnalFunc}. 

Next in section \ref{sec:MTGSim} we show that equivalence of transfer function for two vessels is equivalent to quasi-similarity between them.
After that in section \ref{sec:Adj} the notion of adjoint system is presented. At section \ref{sec:Operations} we present basic operations on vessels:
cascade connection, projection, compression, cascade decomposition. For the completeness of presentation we also discuss in section
\ref{sec:Kalman} Kalman decomposition in our setting that is very similar to the classical $1D$ case.
\section{Overdetermined $t_1$ invariant $2D$ systems}
\subsection{$2D$ systems invariant in one direction}
An overdetermined $t_1$-invariant $2D$ system is a linear input-state-output (i/s/o) system,
which consists of operators depending only on the variable $t_2$; in the most general case such a
system is of the form \cite{bib:Vortices}
\begin{equation*}
  I\Sigma': \left\{ \begin{array}{lll}
         x(t_1, t_2) = e^{A_1(t_2)(t_1-t_1^0)} x(t_1^0, t_2) + \int\limits_{t_1^0}^{t_1}
       e^{A_1(t_2)(t_1-y)} B_1(t_2) u(y, t_2) dy \\
    x(t_1, t_2) = F (t_2,t_2^0) x(t_1, t_2^0) + \int\limits_{t_2^0}^{t_2} F(t_2, s) B_2(s) u(t_1, s)ds \\
    y(t_1, t_2) = C(t_2) x(t_1, t_2) + D(t_2) u(t_1, t_2)
  \end{array} \right.
  \end{equation*}
where for Hilbert spaces $\mathcal E, \mathcal E_*, \mathcal H_{t_2}$ there are defined
\[ \begin{array}{lll}
u(t_1,t_2) \in \mathcal E \text{ - input,} \\
y(t_1,t_2) \in \mathcal E_* \text{ - output,} \\
x(t_1,t_2) \in \mathcal H_{t_2} \text{ - state,}
\end{array}
\]
such that $u(t_1,t_2), y(t_1,t_2)$ are absolutely continuous
functions of each variable when the other variable is fixed. The
transition of the system will usually be considered from
$(t_1^0,t_2^0)$ to $(t_1, t_2)$. Note that $\mathcal H_{t_2}$ are a priory
different for each $t_2$, and as a result $F(t_2,t_2^0)$ has to be
an evolution semi-group, i.e., it satisfies the following definition
\begin{defn} \label{def:EvolSemGr} Given a collection of Hilbert spaces $\{ \mathcal H_t \mid t\in I\}$ for an interval $I\subseteq\mathbb R$
and a collection of bounded invertible operators $F(s,t):\mathcal H_s \rightarrow \mathcal H_t$ for each
$s,t \in I$, we will say that $F(s,t)$ is \textbf{evolution semi-group} if the following relations
hold for all $r,t,s \in I$:
\[ \begin{array}{llll}
F(r,s) F(s,t) = F(r,t), \\
F(t,t) = Id|_{\mathcal H_t}.
\end{array}\]
\end{defn}
In order to be sure that all the formulas are meaningful we shall make the following regularity assumptions
\begin{assump} \label{assm:inner} \textbf{Internal} regularity:
    \begin{enumerate}
    \item $A_1(t_2):\mathcal H_{t_2} \rightarrow \mathcal H_{t_2}$,
       $B_1(t_2), B_2(t_2): \mathcal E \rightarrow \mathcal H_{t_2}$, $ C(t_2): \mathcal H_{t_2} \rightarrow \mathcal E_*$
       are bounded operators (for all $t_2$) and
       $F(t_2,t_2^0): \mathcal H_{t_2^0} \rightarrow \mathcal{H}_{t_2}$ is an evolution semi-group
       (see definition \ref{def:EvolSemGr}).
    \item $F(t_2, s) B_2(s)$ and $C(s) F(s, t_2)$ are absolutely continuous as functions
       of $s$ (for almost all $t_2$) in the norm operator topology on $\mathcal L(\mathcal E, \mathcal H_{t_2})$
       and on $\mathcal L(\mathcal E_*, \mathcal H_{t_2})$, respectively.
    \end{enumerate}
\end{assump}
\begin{assump} \label{assm:FDThrough} \textbf{Feed through} regularity: the operator
$D(t_2): \mathcal E \rightarrow \mathcal E_*$ is an absolutely continuous function of $t_2$.
\end{assump}
Since our entries $u(t_1,t_2)$ will be locally integrable as functions of $t_1$, and as a result the first equation may be equivalently
considered in the differential form, we shall work with systems of the following form:
\begin{equation} \label{eq:systempre}
	I\Sigma: \left\{ \begin{array}{lll}
		\frac{\partial}{\partial t_1}x(t_1,t_2) = A_1(t_2) x(t_1,t_2) + B_1(t_2) u(t_1,t_2) \\
		x(t_1, t_2) = F (t_2,t_2^0) x(t_1, t_2^0) + \int\limits_{t_2^0}^{t_2} F(t_2, s) B_2(s) u(t_1, s)ds \\
		y(t_1, t_2) = C(t_2) x(t_1, t_2) + D(t_2) u(t_1, t_2)
	\end{array} \right.
\end{equation}

\subsection{Overdeterminedness and compatibility}
To ensure that the overdetermined systems equations (\ref{eq:systempre}) are compatible, we shall
demand the equality of the two transitions for our system:
\[
\begin{array}{llllll}
1. ~~~~~ (t_1^0, t_2^0) \longrightarrow (t_1^0, t_2) \longrightarrow (t_1, t_2), \\
2. ~~~~~ (t_1^0, t_2^0) \longrightarrow (t_1, t_2^0) \longrightarrow (t_1, t_2). \\
\end{array}
\]
for arbitrary $(t_1^0,t_2^0), (t_1,t_2)$. In the first case
\begin{multline*}
x(t_1, t_2) = e^{A_1(t_2)(t_1-t_1^0)} F (t_2,t_2^0) x(t_1^0, t_2^0) + \\
       + e^{A_1(t_2)(t_1-t_1^0)} \int\limits_{t_2^0}^{t_2} F(t_2, s) B_2(s) u(t_1^0, s)ds + \\
   + \int\limits_{t_1^0}^{t_1} e^{A_1(t_2)(t_1-p)} B_1(t_2) u(p, t_2) dp
\end{multline*}
and in the second:
\begin{multline*}
\tilde x(t_1, t_2) = F (t_2,t_2^0) e^{A_1(t_2^0)(t_1-t_1^0)} x(t_1^0, t_2^0) + \\
    + F (t_2,t_2^0) \int\limits_{t_1^0}^{t_1} e^{A_1(t_2^0)(t_1-p)} B_1(t_2^0) u(p, t_2^0) dp + \\
  + \int\limits_{t_2^0}^{t_2} F(t_2, s) B_2(s) u(t_1, s)ds.
\end{multline*}
The compatibility condition for free evolution ($u \equiv 0$) results in
\[ e^{A_1(t_2)(t_1-t_1^0)} F (t_2,t_2^0) x(t_1^0, t_2^0) =
F (t_2,t_2^0) e^{A_1(t_2^0)(t_1-t_1^0)} x(t_1^0, t_2^0)
\]
or
\renewcommand{\theequation}{\arabic{section}.\arabic{equation}.Lax}
\begin{equation} \label{eq:LaxCond}
A_1(t_2) = F(t_2, t_2^0) A_1(t_2^0) F(t_2^0, t_2)
\end{equation}
which is called the Lax equation \cite{bib:LaxEquation} and plays an
important role in the theory of completely integrable non-linear
PDEs. Note that it follows that the spectrum of $A_1(t_2)$ is independent of $t_2$.

Inserting the Lax condition into $x(t_1,t_2) = \tilde x(t_1,t_2)$ and rearranging the summands we obtain:
\begin{multline*}
e^{A_1(t_2)(t_1-t_1^0)} \int\limits_{t_2^0}^{t_2} F(t_2, s) B_2(s) u(t_1^0, s)ds - 
\int\limits_{t_2^0}^{t_2} F(t_2, s) B_2(s) u(t_1, s)ds = \\
= F (t_2,t_2^0) \int\limits_{t_1^0}^{t_1} e^{A_1(t_2^0)(t_1-p)} B_1(t_2^0) u(p, t_2^0) dp -
\int\limits_{t_1^0}^{t_1} e^{A_1(t_2)(t_1-p)} B_1(t_2) u(p, t_2) dp
\end{multline*}
Multiplying this equality on the left by $e^{-A_1(t_2)t_1}$ we reach
\begin{multline*}
e^{-A_1(t_2)t_1^0} \int\limits_{t_2^0}^{t_2} F(t_2, s) B_2(s) u(t_1^0, s)ds - 
e^{-A_1(t_2)t_1}\int\limits_{t_2^0}^{t_2} F(t_2, s) B_2(s) u(t_1, s)ds = \\
= e^{-A_1(t_2)t_1}F (t_2,t_2^0) \int\limits_{t_1^0}^{t_1} e^{A_1(t_2^0)(t_1-p)} B_1(t_2^0) u(p, t_2^0) dp -
\int\limits_{t_1^0}^{t_1} e^{-A_1(t_2)p} B_1(t_2) u(p, t_2) dp
\end{multline*}
Since $u(p,s), e^{-A_1(t_2) p}$ are absolutely continuous functions in each variable, we can rewrite it as an equality of iterated integrals
of the derivatives (which are locally absolutely integrable as functions of one variable):
\[ \begin{array}{llllll}
   \int\limits_{t_2^0}^{t_2} \int\limits_{t_1^0}^{t_1}
   \frac{d}{dp} \big[ e^{-A_1(t_2)p} F(t_2, s) B_2(s) u(p, s) \big] dp ds =
   \int\limits_{t_1^0}^{t_1} \int\limits_{t_2^0}^{t_2}
   \frac{d}{ds} \big[ e^{-A_1(t_2) p} F (t_2, s) B_1(s) u(p, s)\big] ds
   dp.
\end{array} \]
Notice that this integral equality is correct for all $t_1^0, t_2^0, t_1, t_2$. Moreover, the functions \linebreak
$e^{-A_1(t_2)p} F(t_2, s) B_2(s) u(p, s)$ and $e^{-A_1(t_2) p} F(t_2, s) B_1(s) u(p, s)$ are absolutely continuous and as a result
their derivatives are integrable functions, meaning that these
two integrals are actually equal by Fubini's theorem to a
two-dimensional integral over the rectangle $\{(p,s) \mid t_1^0 \leq p
\leq t_1, t_2^0 \leq s \leq t_2\}$. Thus it is equivalent to
\begin{gather*} \frac{d}{dp} \big[ e^{-A_1(t_2)p} F(t_2, s) B_2(s) u(p, s) \big] =
   \frac{d}{ds} \big[ e^{-A_1(t_2) p} F (t_2, s) B_1(s) u(p, s)\big]
\end{gather*}
for almost all $(p, s)$. This in turn is simply the following
equation
\begin{gather*} -A_1(t_2) e^{-A_1(t_2) p} F(t_2, s) B_2(s) u(p, s) + e^{-A_1(t_2)p} F(t_2, s) B_2(s)
\frac{\partial}{\partial p}u(p, s) = \\ e^{-A_1(t_2) p} \frac{\partial}{\partial s} [F (t_2, s) B_1(s)] u(p, s) +
e^{-A_1(t_2) p} F (t_2, s) B_1(s) \frac{\partial}{\partial s} u(p, s).
\end{gather*}
Using (\ref{eq:LaxCond}) again we shall obtain
\begin{gather*} - F(t_2, s) A_1(s) B_2(s) u(p, s) + F(t_2, s) B_2(s)
\frac{\partial}{\partial p}u(p, s) = \\
~~~~~~~ = \frac{\partial}{\partial s} [F (t_2, s) B_1(s)] u(p, s)
+ F (t_2, s) B_1(s) \frac{\partial}{\partial s} u(p, s),
\end{gather*}
or, after multiplying on the left by $F(s,t_2)$ and then substituting back the variables $(p,s,t_2) \rightarrow (t_1,t_2,t_2^0)$
\renewcommand{\theequation}{\arabic{section}.\arabic{equation}}
\begin{gather} \label{eq:InpBeforFact}
B_2(t_2) \frac{\partial}{\partial t_1}u(t_1, t_2) -  B_1(t_2) \frac{\partial}{\partial t_2} u(t_1, t_2)
- \big( A_1(t_2) B_2(t_2) + F(t_2, t_2^0)\frac{\partial}{\partial t_2} [F (t_2^0, t_2) B_1(t_2)] \big)
u(t_1, t_2) = 0.
\end{gather}

At this stage it is convenient to assume that we have factorization
\begin{equation} \label{eq:Factorization}
\begin{array}{lllll}
B_2(t_2) = \widetilde B(t_2) \sigma_2(t_2), ~~B_1(t_2) = \widetilde B(t_2) \sigma_1(t_2), \\
A_1(t_2) B_2(t_2) + F(t_2, t_2^0)\frac{\partial}{\partial s} [F (t_2^0, t_2) B_1(t_2)] = - \widetilde B(t_2) \gamma(t_2)
\end{array} \end{equation}
for some operators
\[ \widetilde B(t_2): \widetilde{\mathcal E} \rightarrow \mathcal H_{t_2}, ~~
\sigma_2(t_2), \sigma_1(t_2), \gamma(t_2): \mathcal E \rightarrow \widetilde{\mathcal E},
\]
where $ \widetilde{\mathcal E}$ is another auxiliary Hilbert
space. It is also important to postulate the following assumption in
order to give meaning to the corresponding formulas:
\begin{assump} \label{assm:extin} \textbf{External input} regularity:
        \begin{enumerate}
        \item $\gamma(t_2), \sigma_2(t_2) \in L^1_{loc}( \mathcal L(\mathcal E, \widetilde{\mathcal E}))$
                in the norm operator topology.
        \item $\sigma_1(t_2) \in L(\mathcal E, \widetilde{\mathcal E})$ is absolutely continuous and invertible,
                        in the norm operator topology.
        \end{enumerate}
\end{assump}

Expressed directly in terms of the operators $\widetilde B(t_2),
\sigma_1(t_2), \sigma_2(t_2)$ with $B_1(t_2), B_2(t_2)$
eliminated, (\ref{eq:Factorization}) becomes
\renewcommand{\theequation}{\arabic{section}.\arabic{equation}.OverD}
\begin{gather} \label{eq:OverDetCondIn}
\frac{d}{dt_2} (F (t_2^0, t_2) \widetilde B(t_2) \sigma_1(t_2)) + F(t_2^0, t_2) A_1(t_2) \widetilde B(t_2) \sigma_2(t_2)
   + F(t_2^0, t_2) \widetilde B(t_2) \gamma(t_2) = 0.
\end{gather}
Then the condition (\ref{eq:InpBeforFact}) becomes
\[ \widetilde B(t_2) [\sigma_2(t_2) \frac{\partial}{\partial t_1}u(t_1, t_2) -
  \sigma_1(t_2) \frac{\partial}{\partial t_2}u(t_1,t_2) + \gamma(t_2)] u(t_1,t_2) = 0.
\]
A sufficient condition for this to hold (which is necessary in case $\widetilde B(t_2)$ is injective) is
the \textit{input compatibility condition}
\renewcommand{\theequation}{\arabic{section}.\arabic{equation}}
\begin{equation}\label{eq:InCCP}
\sigma_2(t_2) \frac{\partial}{\partial t_1}u(t_1, t_2) -
  \sigma_1(t_2) \frac{\partial}{\partial t_2}u(t_1,t_2) + \gamma(t_2) u(t_1,t_2) = 0
\end{equation}
Note that since $u(t_1,t_2)$ is the solution of PDE in the extended sense, it will be absolutely
continuous as a function of $t_1$ for almost all $t_2$, and conversely, it will be absolutely
continuous as a function of $t_2$ for almost all $t_1$.

The output $y(t_1,t_2)$ should satisfy the \textit{output compatibility condition} of the same type as for the input compatibility
condition (\ref{eq:InCCP}), namely:
\begin{equation} \label{eq:OutCCP}
\sigma_{2*}(t_2) \frac{\partial}{\partial t_1}y(t_1, t_2) -
  \sigma_{1*}(t_2) \frac{\partial}{\partial t_2}y(t_1,t_2) + \gamma_*(t_2) y(t_1,t_2) = 0.
\end{equation}
where similarly we have the following assumptions using an auxiliary output space $\widetilde{\mathcal E_*}$
\begin{assump} \label{assm:extout} \textbf{External output} regularity:
        \begin{enumerate}
        \item $\gamma_*(t_2), \sigma_{2*}(t_2) \in L^1_{loc}( \mathcal L(\mathcal E_*, \widetilde{\mathcal E_*}))$
                in the norm operator topology.
        \item $\sigma_{1*}(t_2) \in L(\mathcal E_*, \widetilde{\mathcal E_*})$ is absolutely continuous and invertible,
                        in the norm operator topology.
        \end{enumerate}
\end{assump}

So, inserting here $y(t_1, t_2) = D(t_2)u(t_1,t_2) + C(t_2) x(t_2,t_2)$ we obtain that
\[ \begin{array}{lllllllll}
0 & = [\sigma_{2*}(t_2) \frac{\partial}{\partial t_1} - \sigma_{1*}(t_2) \frac{\partial}{\partial t_2} + \gamma_*(t_2)] y(t_1,t_2) = \\
    & = [\sigma_{2*}(t_2) \frac{\partial}{\partial t_1} - \sigma_{1*}(t_2) \frac{\partial}{\partial t_2} + \gamma_*(t_2)][D(t_2)u(t_1,t_2) + C(t_2) x(t_2,t_2)]\\
    & = \sigma_{2*}(t_2) C(t_2) \frac{\partial}{\partial t_1}x(t_2,t_2) - \sigma_{1*}(t_2)\frac{\partial}{\partial t_2}[C(t_2)x(t_2,t_2)] + \gamma_*(t_2)
    C(t_2) x(t_2,t_2) + \\
    & ~~~~~~ + \sigma_{2*}(t_2) D(t_2) \frac{\partial}{\partial t_1}u(t_2,t_2) - \sigma_{1*}(t_2) \frac{\partial}{\partial t_2} [D(t_2) u(t_2,t_2)] + \gamma_*(t_2) D(t_2) u(t_2,t_2).
\end{array} \]
Substituting here the first system equation from
(\ref{eq:systempre}), we obtain (after omitting the notation
of dependence on the variables)
\[ \begin{array}{lllllllll}
0 & = \sigma_{2*} C (A_1 x + \widetilde B \sigma_1 u) - \sigma_{1*}\frac{\partial}{\partial t_2}[C x] + \gamma_*    C x + \sigma_{2*} D
        \frac{\partial}{\partial t_1}u - \sigma_{1*} [D' u + D \frac{\partial}{\partial t_2} u ] + \gamma_*D u = \\
    & = \sigma_{2*} C A_1 x - \sigma_{1*}\frac{\partial}{\partial t_2}[C x] + \gamma_*  C x + \\
    & ~~~~~~~~~~~~ +    [\sigma_{2*} D \frac{\partial}{\partial t_1} - \sigma_{1*} D \frac{\partial}{\partial t_2} +
            \sigma_{2*} C \widetilde B \sigma_1 - \sigma_{1*} D' + \gamma_*D] u.
\end{array} \]
Multiplying the second equation of the system (\ref{eq:systempre}) by $C(t_2)$ and differentiating with respect to $t_2$,
one can easily obtain that $\frac{\partial}{\partial t_2}[C x]$ satisfies:
\[
\frac{\partial}{\partial t_2}[C x] = \frac{\partial}{\partial
t_2}[C F] F^{-1}x + C \widetilde B \sigma_2 u.
\]
Inserting this into the last equation we obtain
\begin{equation} \label{eq:OutGen} \begin{array}{lllllllll}
0 & =  \sigma_{2*} C A_1 x - \sigma_{1*}[ \frac{\partial}{\partial t_2}[C F] F^{-1}x + C \widetilde B \sigma_2 u] + \gamma_*    C x + \\
    & ~~~~~~~~~~~~ +    [\sigma_{2*} D \frac{\partial}{\partial t_1} - \sigma_{1*} D \frac{\partial}{\partial t_2} +
            \sigma_{2*} C \widetilde B \sigma_1  - \sigma_{1*} D' + \gamma_*D] u = \\
    & = [ \sigma_{2*} C A_1 - \sigma_{1*} \frac{\partial}{\partial t_2}[C F] F^{-1} + \gamma_*  C] x  + \\
    & ~~~~~~~~~~~~ +    [\sigma_{2*} D \frac{\partial}{\partial t_1} - \sigma_{1*} D \frac{\partial}{\partial t_2} +
            \sigma_{2*} C \widetilde B \sigma_1  - \sigma_{1*} C \widetilde B \sigma_2 + \sigma_{1*} D' + \gamma_*D] u.
\end{array} \end{equation}
The validity of this equation for the special case $u(t_1,t_2) = 0$ and an arbitrary initial $x(t_1^0,t_2^0)$ forces us to impose
\renewcommand{\theequation}{\arabic{section}.\arabic{equation}.OverD}
\begin{gather} \label{eq:OverDetCondOut}
0 = \sigma_{2*} C A_1 F  - \sigma_{1*} \frac{\partial}{\partial t_2}[C F] + \gamma_* C F.
\end{gather}
With (\ref{eq:OverDetCondOut}) in force,  (\ref{eq:OutGen}) collapses to
\[ 0 = [\sigma_{2*} D \frac{\partial}{\partial t_1} - \sigma_{1*} D \frac{\partial}{\partial t_2} +
            \sigma_{2*} C \widetilde B \sigma_1  - \sigma_{1*} C \widetilde B \sigma_2 + \sigma_{1*} D' + \gamma_*D] u.
\]
On the other hand, $u$ satisfies the input compatibility condition
(\ref{eq:InCCP}), so it is natural to assume that there is an
operator $\widetilde D: \mathcal E_* \rightarrow \mathcal
{\widetilde E_*}$ satisfying
\begin{assump} \label{assm:TildeD} \textbf{External feed through} regularity:  the operator
$\widetilde D(t_2):  \mathcal E_* \rightarrow \mathcal{\widetilde E_*}$ is an absolutely continuous function of $t_2$,
\end{assump}
and which satisfies the following intertwining conditions, which will be called from now on the \textit{linkage
conditions}
\renewcommand{\theequation}{\arabic{section}.\arabic{equation}.Link}
\begin{equation} \label{eq:LinkCond}
\begin{array}{ll}
\sigma_{1*} D = \widetilde D \sigma_1, ~~~~~~ \sigma_{2*} D = \widetilde D \sigma_2,  \\
\widetilde D \gamma = \sigma_{2*} C \widetilde B \sigma_1 - \sigma_{1*} C \widetilde B \sigma_2 + \sigma_{1*} D' + \gamma_*D.
\end{array} \end{equation}

\section{\label{sec:Vessel}Vessel}
\subsection{Definition}
Let us combine together all the formulas we have just developed.
We define an (integral) \textbf{vessel} to be a collection of operator and
spaces
\[ \mathfrak{IV} = (A_1(t_2), F(t_2,t_2^0),\widetilde B(t_2),C(t_2), D(t_2),\widetilde D(t_2); \sigma_1(t_2), \sigma_2(t_2), \gamma(t_2),
\sigma_{1*}(t_2), \sigma_{2*}(t_2)\gamma_*(t_2); \mathcal H_{t_2}, \mathcal{E},\mathcal E_*, \mathcal{\widetilde E},\mathcal{\widetilde E_*})
\]
satisfying regularity assumptions \ref{assm:inner},
\ref{assm:FDThrough}, \ref{assm:extin}, \ref{assm:extout}, \ref{assm:TildeD} and the following vessel
conditions:
\[ \begin{array}{lllllll}
 F(t_2,t_2^0) A_1(t_2^0) = A_1(t_2) F(t_2,t_2^0)  & (\text{\ref{eq:LaxCond}}) \\
  \frac{d}{dt_2} (F (t_2^0, t_2) \widetilde B(t_2) \sigma_1(t_2)) + F(t_2^0, t_2) A_1(t_2) \widetilde B(t_2) \sigma_2(t_2)
   + F(t_2^0, t_2) \widetilde B(t_2) \gamma(t_2) = 0
              & (\text{\ref{eq:OverDetCondIn}}) \\
  \sigma_{2*}(t_2) C(t_2) A_1(t_2) F(t_2,t_2^0)  - \sigma_{1*}(t_2) \frac{d}{d t_2}[C(t_2) F(t_2,t_2^0)] +
        \gamma_*(t_2) C(t_2) F(t_2,t_2^0) = 0             & (\text{\ref{eq:OverDetCondOut}}) \\
 \begin{array}{ll}
\sigma_{1*} D = \widetilde D \sigma_1, ~~~~~~ \sigma_{2*} D = \widetilde D \sigma_2,  \\
\widetilde D \gamma = \sigma_{2*} C \widetilde B \sigma_1  - \sigma_{1*} C \widetilde B \sigma_2 + \sigma_{1*} D' + \gamma_*D.
\end{array}   & (\text{\ref{eq:LinkCond}})
\end{array} \]
It is naturally associated to the system $I\Sigma$ (see (\ref{eq:systempre}))
\renewcommand{\theequation}{\arabic{section}.\arabic{equation}}
\begin{equation} \label{eq:system}
     I\Sigma: \left\{ \begin{array}{lll}
    \frac{\partial}{\partial t_1}x(t_1,t_2) = A_1(t_2) ~x(t_1,t_2) + \widetilde B(t_2) \sigma_1(t_2) ~u(t_1,t_2) \\[5pt]
    x(t_1, t_2) = F (t_2,t_2^0) x(t_1, t_2^0) + \int\limits_{t_2^0}^{t_2} F(t_2, s) \widetilde B(s) \sigma_2(s) u(t_1, s)ds \\[5pt]
    y(t_1,t_2) = C(t_2)~ x(t_1,t_2) + D(t_2) u(t_1,t_2).
  \end{array} \right.
  \end{equation}
with absolutely continuous inputs and outputs, satisfying
compatibility conditions (\ref{eq:InCCP}), (\ref{eq:OutCCP}) for
almost all $(t_1,t_2)$:
\[ \begin{array}{llll}
  \sigma_2(t_2) \frac{\partial}{\partial t_1}u(t_1, t_2) -
  \sigma_1(t_2) \frac{\partial}{\partial t_2}u(t_1,t_2) + \gamma(t_2) u(t_1,t_2) = 0, \\
  \sigma_{2*}(t_2) \frac{\partial}{\partial t_1}y(t_1, t_2) -
  \sigma_{1*}(t_2) \frac{\partial}{\partial t_2}y(t_1,t_2) + \gamma_*(t_2) y(t_1,t_2) = 0.
\end{array} \]
We shall further name these conditions as follows.
(\ref{eq:LaxCond}) is the \textit{Lax equation}.
(\ref{eq:OverDetCondIn}) and (\ref{eq:OverDetCondOut}) are the
\textit{input and the output vessel conditions}, respectively.
(\ref{eq:LinkCond}) is the \textit{linkage condition}.

\subsection{Gauge quasi-similarity of vessels}
As in the classical case, in order to deal with some classification
of systems (to be defined later) we need a notion of minimal
systems. In the $1D$ case there is only one natural notion of
approximate controllability and observability. In the case of
$2D$ $t_1$ invariant systems, on the other hand, there are at least the following two notions
\begin{defn} System $I\Sigma$ (\ref{eq:system}) is called \textbf{locally approximately controllable at $t_2$} if
\[ \mathcal C_{t_2} = \{ h\in\mathcal H_{t_2} \mid \exists t_1\in \mathbb R, (u,x,y)\in\mathcal T: x(0,t_2) = 0, x(t_1,t_2) = h  \}
\]
is dense in $\mathcal H_{t_2}$.
\end{defn}
Here $\mathcal{T}$ stands for the set of system trajectories $(u,x,y)$ of (\ref{eq:system}) with compatibility
ODEs \ref{eq:InCC}) and (\ref{eq:OutCC}) with the spectral parameter $\lambda$.
\begin{defn} System $I\Sigma$ (\ref{eq:system}) is called \textbf{approximately controllable at $t_2$} if
\[ \widetilde{\mathcal C}_{t_2} = \{ h\in\mathcal H_{t_2} \mid \exists t_1\in \mathbb R, (u,x,y)\in\mathcal T: x(0,0) = 0, x(t_1,t_2) = h  \}
\]
is dense in $\mathcal H_{t_2}$.
\end{defn}
Similarly to the classical case one immediately obtains that
\[ \begin{array}{ll}
 \mathcal C_{t_2} = \operatorname{\overline{span}} \{ \operatorname{Im} A_1^j(t_2) B(t_2) \mid j=0,1,2, \ldots \}, \\
 \widetilde{\mathcal C}_{t_2} = \operatorname{\overline{span}} \{ \operatorname{Im} F(t_2,s) A_1^j(s) B(s) \mid j=0,1,2, \ldots, s\in\mathbb R \}. \\
\end{array} \]
Notice that $F(t_2,t_2') \widetilde{\mathcal C}_{t_2'} = \widetilde{\mathcal C}_{t_2}$ and since
$F(t_2,t_2')$ are bounded invertible operators, the density of $\widetilde{\mathcal C}_{t_2}$ for any value of $t_2$
implies density for all $t_2$, which means that the notion of approximate controllability is independent of $t_2$.

Since $F(t_2,t_2) = Id$, we obtain that $\mathcal C_{t_2} \subseteq \widetilde{\mathcal C}_{t_2}$ for all $t_2$
and consequently, if $\mathcal C_{t_2}$ is dense in $\mathcal H_{t_2}$ then so is $\widetilde{\mathcal C}_{t_2}$.
But actually the converse also holds
\begin{thm} \label{thm:GlobLocContrl}
For the system $I\Sigma$ (\ref{eq:system}), satisfying regularity assumptions \ref{assm:inner},
\ref{assm:FDThrough}, \ref{assm:extin}, \ref{assm:extout}, \ref{assm:TildeD}
the local approximate controllability for arbitrary $t_2^0$ is equivalent to approximate controllability for all $t_2$.
\end{thm}
\textbf{Proof:} We have seen that
\[ \mathcal C_{t_2^0} = \mathcal H_{t_2^0} \Rightarrow  \widetilde{\mathcal C}_{t_2^0} = \mathcal H_{t_2^0} \Rightarrow
\widetilde{\mathcal C}_{t_2} = \mathcal H_{t_2}
\]
and we need to prove the converse. Suppose that for a fixed $t_2^0$ we have approximate controllability, which means that
\[ \bigvee_{n\geq 0,e \in \mathcal E} F(t_2^0,s) A_1^j(s) B(s) e = \mathcal H_{t_2^0}
\]
or what is equivalent
\[ B^*(s) F^*(t_2^0,s) A_1^*(t_2^0) h = 0, ~~\forall h\in\mathcal H_{t_2^0}.
\]
In the same manner local approximate controllability means that
\[ \bigvee_{n\geq 0,e \in \mathcal E} A_1^n(t_2^0) B(t_2^0) e = \mathcal H_{t_2^0} \Rightarrow
B^*(t_2^0) A_1^{*n}(t_2^0) h = 0, ~~\forall h\in\mathcal H_{t_2^0}.
\]
Consider now a function with values in $\mathcal E$ for each $h\in\mathcal H_{t_2^0}$
\[ u_h(\lambda,t_2,t_2^0) = B^*(t_2) F^*(t_2^0,t_2) (\lambda I - A_1^*(t_2^0))^{-1}h, ~~u_h(\lambda,t_2^0,t_2^0) = B^*(t_2^0) (\lambda I - A_1^*(t_2^0))^{-1}h
\]
which is analytic at the neighborhood of $\lambda = \infty$. This function satisfies the following differential equation 
(adjoint of the input vessel condition (\ref{eq:OverDetCondIn}))
\[ \begin{array}{lll}
\frac{d}{dt_2} u_h(\lambda,t_2,t_2^0) = \frac{d}{dt_2} [B^*(t_2) F^*(t_2^0,t_2) (\lambda I - A_1^*(t_2^0))^{-1}h] = \\
~~~~~= \sigma_1^{-1}[-\sigma_2 B^*(t_2) F^*(t_2^0,t_2) A_1^{*n}(t_2^0) (\lambda I - A_1^*(t_2^0))^{-1}h +
				\gamma B^*(t_2) F^*(t_2^0,t_2) A_1^{*n}(t_2^0) (\lambda I - A_1^*(t_2^0))^{-1}h] = \\
~~~~~= \sigma_1^{-1}[-\sigma_2 B^*(t_2) F^*(t_2^0,t_2) (A_1^{*n}(t_2^0) -\lambda I + \lambda I)(\lambda I - A_1^*(t_2^0))^{-1}h +
				\gamma u_h(\lambda,t_2,t_2^0)] = \\
~~~~~= \sigma_1^{-1}[-\sigma_2 \lambda + \gamma] u_h(\lambda,t_2,t_2^0) +	\sigma_1^{-1}\sigma_2 B^*(t_2) F^*(t_2^0,t_2) h.
\end{array} \]
denoting by $\Phi(\lambda,t_2,t_2^0)$ the fundamental solution of the differential equation
\[ \frac{d}{dt_2} \Phi(\lambda,t_2,t_2^0) = \sigma_1^{-1}[-\sigma_2 \lambda + \gamma] \Phi(\lambda,t_2,t_2^0)
\]
and using variation of parameters, we shall obtain that
\[ \begin{array}{lll}
u_h(\lambda,t_2,t_2^0) & = \Phi(\lambda,t_2,t_2^0)[\int_{t_2^0}^{t_2}\Phi^{-1}(\lambda,y,t_2^0)\sigma_1^{-1}\sigma_2 B^*(y) F^*(t_2^0,y) h +
u_h(\lambda,t_2^0,t_2^0)] = \\
& = \Phi(\lambda,t_2,t_2^0)[\int_{t_2^0}^{t_2}\Phi^{-1}(\lambda,y,t_2^0)\sigma_1^{-1}\sigma_2 B^*(y) F^*(t_2^0,y) h +
B^*(t_2^0) (\lambda I - A_1^n(t_2^0))^{-1}h] = \\
& = \Phi(\lambda,t_2,t_2^0)\int_{t_2^0}^{t_2}\Phi^{-1}(\lambda,y,t_2^0)\sigma_1^{-1}\sigma_2 B^*(y) F^*(t_2^0,y) h +
\Phi(\lambda,t_2,t_2^0)B^*(t_2^0) (\lambda I - A_1^n(t_2^0))^{-1}h,
\end{array} \]
which is a sum of an analytic (in $\lambda$) $\mathcal E$-valued function and of a function with poles (for $\lambda$ in the spectrum of $A_1^*(t_2^0)$)
\[ \Phi(\lambda,t_2,t_2^0)B^*(t_2^0) (\lambda I - A_1^*(t_2^0))^{-1}h.
\]
So, if there is a vector $h\in\mathcal H_{t_2^0}$ such that $u_h(\lambda,t_2,t_2^0) = 0$ for each value of $\lambda$, it has to vanish
the part containing the pole first. This means that for all $\lambda$ out of the spectrum of $A_1^*(t_2^0)$
\[ \Phi(\lambda,t_2,t_2^0)B^*(t_2^0) (\lambda I - A_1^*(t_2^0))^{-1}h = 0 
\]
and since $\Phi(\lambda,t_2,t_2^0)$ is invertible we obtain,
\[ B^*(t_2^0) A_1^{*n}(t_2^0) h = 0,
\]
which is local approximate controllability.
\qed

The following notions are natural generalizations of the $1D$ case:
\begin{defn} The system $I\Sigma$ is \textbf{locally observable at $t_2$} if
\[ \mathcal O^\perp_{t_2} = \{ h \in \mathcal H_{t_2} \mid \text{$x(0,t_2) = h$,
$(u,x,y) \in \mathcal T$, $u(t_1,t_2) = 0 \Rightarrow y(t_1,t_2) = 0$ for all $t_1$} \} = \{ 0 \},
\]
and the system is \textbf{observable at $t_2$} if
\[ \widetilde{\mathcal O}^\perp_{t_2} = \{ h \in \mathcal H \mid \text{$x(0,0) = h$,
$(u,x,y) \in \mathcal T$, $u = 0 \Rightarrow y = 0$ for all $t_1$} \} = \{ 0 \}.
\]
\end{defn}
Again, similarly to the classical case one obtains that
\[ \begin{array}{lll}
\mathcal O^\perp_{t_2}             &  = \bigcap_{n\in\mathbb N} \operatorname{Ker} C(t_2) A_1^n(t_2), \\
\widetilde{\mathcal O}^\perp_{t_2} & = \bigcap_{n\in\mathbb N, s\in\mathbb R} \operatorname{Ker} C(s) A_1^n(s) F(s,t_2) .
\end{array}\]
From the property of the evolution semi-group $F(t_2,t_2) = I$ we obtain that
$\mathcal O^\perp_{t_2} \supseteq \widetilde{\mathcal O}^\perp_{t_2}$ and that local
observability at $t_2$ (i.e., $\mathcal O^\perp_{t_2} = \{ 0 \})$ implies observability at $t_2$
(i.e., $\widetilde{\mathcal O}^\perp_{t_2} = \{ 0] \}$). Additionally,
$\widetilde{\mathcal O}^\perp_{t_2} = \widetilde{\mathcal O}^\perp_{t_2'} F(t_2',t_2)$ implies that
observability is independent of $t_2$, i.e., if for any $t_2^0$, $\widetilde{\mathcal O}^\perp_{t_2^0} = \{0\}$ then
$\widetilde{\mathcal O}^\perp_{t_2} = \{ 0 \}$ for all $t_2$. Finally, a parallel to the controllability theorem.
\begin{thm} For the system $I\Sigma$ (\ref{eq:system}), satisfying regularity assumptions \ref{assm:inner},
\ref{assm:FDThrough}, \ref{assm:extin}, \ref{assm:extout}, \ref{assm:TildeD} the local observability for arbitrary $t_2^0$ is equivalent to
observability for all $t_2$.
\end{thm}
\textbf{Proof:} We have seen that
\[ \mathcal O^\perp_{t_2^0} = \{ 0\} \Rightarrow  \widetilde{\mathcal O}^\perp_{t_2^0} = \{ 0\}  \Rightarrow
\widetilde{\mathcal O}^\perp_{t_2} = \{ 0\}
\]
and it is remained to show the converse. Similarly to the proof of local approximate controllability for each $t_2$, supposing 
observability at a fixed $t_2$
\[ \widetilde{\mathcal O}^\perp_{t_2} = \bigcap_{n\in\mathbb N, s\in\mathbb R} \operatorname{Ker} C(s) A_1^n(s) F(s,t_2) = \{ 0 \},
\]
we may consider a function of $\lambda$ with values in $\mathcal E$
\[ y_h(\lambda,t_2,t_2^0) = C(t_2) F(t_2,t_2^0) (\lambda I - A_1(t_2^0))^{-1}h,
\]
and the output vessel condition (\ref{eq:OverDetCondOut})
\[ \sigma_{2*}(t_2) y_{A_1(t_2^0)h}(t_2)  - \sigma_{1*}(t_2) \frac{d}{d t_2}y_h(t_2) + \gamma_*(t_2) y_h(t_2) = 0, 
~~y_h(\lambda,t_2^0,t_2^0) = C(t_2^0) (\lambda I - A_1(t_2^0))^{-1}h
\]
and follow the same lines as in the proof of local approximate controllability.   \qed

The system $I\Sigma$ is called \textit{minimal} if it is both approximately controllable and observable, i.e.,
($\mathcal{C}_{t_2} = \mathcal{O}_{t_2} = \mathcal{H}$).

A natural notion of similarity arises, which is used to classify (usually minimal) systems and corresponding vessels. Two vessels
\[ \begin{array}{lll}
 \mathfrak{IV} = (A_1(t_2), F(t_2,t_2^0),\widetilde B(t_2),C(t_2), D(t_2),\widetilde D(t_2); \sigma_1(t_2), \sigma_2(t_2), \gamma(t_2),
\sigma_{1*}(t_2), \sigma_{2*}(t_2)\gamma_*(t_2); \mathcal H_{t_2}, \mathcal{E},\mathcal E_*, \mathcal{\widetilde E},\mathcal{\widetilde E_*})
 \\
 \mathfrak{\breve{IV}} = (\breve{A}_1(t_2), \breve F(t_2,t_2^0),\widetilde{\breve B}(t_2),\breve C(t_2), D(t_2),\widetilde{D}(t_2); \sigma_1(t_2), \sigma_2(t_2), \gamma(t_2),
\sigma_{1*}(t_2), \sigma_{2*}(t_2)\gamma_*(t_2); \breve{\mathcal H}_{t_2}, \mathcal{E},\mathcal E_*, \mathcal{\widetilde E},\mathcal{\widetilde E_*})
\end{array} \]
are called \textit{gauge quasi-similar} (or \textit{kinematically quasi-similar}),
if there exists a (possibly unbounded) linear operator $T(t_2): D(T(t_2)) \rightarrow \breve{\mathcal{H}}_{t_2}$, with
a dense domain $D(T(t_2)) \subseteq \mathcal H_{t_2}$, which is 1-1, with dense range, and satisfies the following intertwining conditions
\begin{equation} \label{eq:Uconnect} \left\{\begin{array}{llll}
\breve A_1(t_2) T(t_2)       & = T(t_2) A_{1}(t_2), \\
\breve F(t_2,t_2^0) T(t_2^0) & = T(t_2) F(t_2,t_2^0), \\
\widetilde {\breve B}(t_2)   & = T (t_2) \widetilde B(t_2), \\
\breve C(t_2) T(t_2)         & = C(t_2).
\end{array}\right.
\end{equation}
Moreover, we shall demand that
\begin{equation} \label{eq:TDomainCond}
\widetilde{\mathcal C}_{t_2} \subseteq D(T(t_2)), ~~~~~~~~ F(t_2,t_2^0) D(T(t_2^0)) \subseteq D(T(t_2)).
\end{equation}
These conditions are necessary in order to obtain reasonable definitions in (\ref{eq:Uconnect}). For example, 
$\breve F(t_2,t_2^0) T(t_2^0) = T(t_2) F(t_2,t_2^0)$ requires $F(t_2,t_2^0) D(T(t_2^0)) \subseteq D(T(t_2))$. The first and the third conditions
require \linebreak $\operatorname{Im}B(t_2), \operatorname{Im} A_1(t_2) \subseteq D(T_{t_2})$, for which it is enough to demand 
$\widetilde{\mathcal C}_{t_2} \subseteq D(T(t_2))$.

When it is the case that $D(T(t_2)) = \mathcal H_{t_2}$ for all $t_2$, and $T(t_2)$ is everywhere defined bounded and onto
(with bounded inverse as a result), we say that the vessels $\mathfrak{IV}, \mathfrak{\breve{IV}}$ are \textit{similar}.

\noindent\textbf{Remark: 1.} For the finite dimensional case $\dim\mathcal E < \infty$, as in the classical case the notions of similarity and of 
quasi-similarity coincide.

\noindent\textbf{2.} Given a vessel $\mathfrak{IV}$ and a family of invertible bounded operators 
$T(t_2) \colon \mathcal H_{t_2} \to \breve{\mathcal H_{t_2}}$ the formulas (\ref{eq:Uconnect}) define a new vessel 
$\mathfrak{\breve{IV}}$ that is gauge similar to $\mathfrak{IV}$.

\subsection{\label{sec:DifForm}Differential vessels}
Suppose that the system $I\Sigma$ defined by (\ref{eq:system}) has identical inner spaces $\mathcal H_{t_2}$ for all $t_2$, i.e.,
$\mathcal H_{t_2} = \mathcal H$. Suppose also that the evolution semi-group $F(t_2,t_2^0)$ is absolutely
continuous (for a fixed $t_2^0$) in the norm operator topology of $B(\mathcal H, \mathcal H)$, then
its generator $A_2(t_2)$ is for almost all $t_2$ a bounded operator and satisfies
\[ A_2(t_2) = \frac{d}{dt_2} F(t_2,t_2^0) F^{-1}(t_2,t_2^0).
\]
Rewriting all the vessel conditions in the differential form using $A_2(t_2)$, we obtain a \textit{differential vessel}
\[ \mathfrak{DV} = (A_1(t_2), A_2(t_2),\widetilde{B}(t_2),C(t_2), D(t_2),\widetilde D(t_2);
        \sigma_1(t_2), \sigma_2(t_2), \gamma(t_2), \sigma_{1*}(t_2), \sigma_{2*}(t_2)\gamma_*(t_2);
        \mathcal H, \mathcal{E},\mathcal E_*, \mathcal{\widetilde E},\mathcal{\widetilde E_*})
\]
which satisfy the following axioms:
\[ \begin{array}{lllllll}
    \frac{d}{dt_2} A_1(t_2) = A_2(t_2) A_1(t_2) - A_1(t_2) A_2(t_2) \\
    \frac{d}{dt_2} \big(\widetilde{B}(t_2) \sigma_1(t_2)\big) - A_2(t_2) \widetilde{B}(t_2) \sigma_1(t_2) + A_1(t_2)
            \widetilde{B}(t_2) \sigma_2(t_2) + \widetilde{B}(t_2) \gamma(t_2) = 0 \\
        \frac{d}{dt_2} \big( \sigma_{1*}(t_2) C(t_2) \big) + \sigma_{1*}(t_2) C(t_2) A_2(t_2) +
                \sigma_{2*}(t_2) C(t_2)A_1(t_2) + \gamma_*(t_2) C(t_2) = 0 \\
 \begin{array}{ll}
\sigma_{1*} D = \widetilde D \sigma_1, ~~~~~~ \sigma_{2*} D = \widetilde D \sigma_2,  \\
\widetilde D \gamma = \sigma_{2*} C \widetilde B \sigma_1  - \sigma_{1*} C \widetilde B \sigma_2 + \sigma_{1*} D' + \gamma_*D.
\end{array}
\end{array} \]
and the following regularity assumptions, which are obtained from the assumptions \ref{assm:inner},
\ref{assm:FDThrough}, \ref{assm:extin}, \ref{assm:extout}, \ref{assm:TildeD} by ,,differentiating''
\begin{assump} 
    \begin{enumerate}
    \item \textbf{Internal} regularity: $A_1(t_2),A_2(t_2):\mathcal H_{t_2} \rightarrow \mathcal H_{t_2}$,
       $B_1(t_2), B_2(t_2): \mathcal E \rightarrow \mathcal H_{t_2}$, $ C(t_2): \mathcal H_{t_2} \rightarrow \mathcal E_*$
       are bounded operators (for all $t_2$),
		\item \textbf{Feed through} regularity: the operator $D(t_2): \mathcal E \rightarrow \mathcal E_*$,
				$\widetilde D(t_2):  \mathcal E_* \rightarrow \mathcal{\widetilde E_*}$ are  
				absolutely continuous function of $t_2$,
		\item \textbf{External regularity} \ref{assm:extin},\ref{assm:extout},
		\end{enumerate}
\end{assump}

The differential vessel is associated with the system
\begin{equation} \label{eq:Dsystem}
    D\Sigma: \left\{ \begin{array}{lll}
    \frac{\partial}{\partial t_1}x(t_1,t_2) = A_1(t_2) x(t_1,t_2) + \widetilde{B}(t_2) ~\sigma_1(t_2) ~u(t_1,t_2) \\
    \frac{\partial}{\partial t_2}x(t_1,t_2) = A_2(t_2) x(t_1,t_2) + \widetilde{B}(t_2) ~\sigma_2(t_2) ~u(t_1,t_2) \\
    y(t_1,t_2) = D(t_2) u(t_1,t_2) + C(t_2) x(t_1,t_2)
    \end{array} \right.
\end{equation}
and compatibility conditions for the input/output signals:
\[ \begin{array}{lll}
\sigma_2(t_2) \frac{\partial}{\partial t_1}u(t_1,t_2) - \sigma_1(t_2) \frac{\partial}{\partial t_2}u(t_1,t_2) +
\gamma(t_2) u(t_1,t_2) = 0, \\
\sigma_{2*}(t_2) \frac{\partial}{\partial t_1}y(t_1,t_2) - \sigma_{1*}(t_2) \frac{\partial}{\partial t_2}y(t_1,t_2) +
\gamma_*(t_2) y(t_1,t_2) = 0.
\end{array} \]
This motivates the following definition
\begin{defn} Vessel $\mathfrak{IV}$ will be called \textbf{differentiable} if $\mathcal H_{t_2} = \mathcal H$ for all $t_2$ and
$F(t_2,t_2^0)$ is an absolutely continuous function in the norm operator topology of $B(\mathcal H, \mathcal H)$.
\end{defn}
Then the following proposition holds
\begin{prop} Any vessel $\mathfrak{IV}$ is gauge similar to a differentiable vessel.
\end{prop}
\textbf{Proof:} Let us take $T(t_2) = F(t_2^0,t_2): \mathcal H_{t_2} \rightarrow \mathcal H_{t_2^0}$, which is an absolutely continous,
and let us find out what kind of a vessel is obtained.
First notice that $T(t_2)$ are all mappings onto the same space $\mathcal H_{t_2^0}$, which we denote by
$\mathcal H_{t_2^0} = \mathcal H$. Further, we see that
\[ \left\{\begin{array}{llll}
\breve A_1(t_2) F(t_2^0,t_2) = F(t_2^0,t_2^0) A_1(t_2)      & \Rightarrow & \breve A_1(t_2) = F(t_2^0,t_2) A_1(t_2) F(t_2,t_2^0), \\
\breve F(t_2,t_2^0) F(t_2^0,t_2^0) = F(t_2^0,t_2) F(t_2,t_2^0)              & \Rightarrow & \breve F(t_2,t_2^0) = I,\\
\widetilde {\breve B}(t_2)   = F(t_2^0,t_2) \widetilde B(t_2),   \\
\breve C(t_2) F(t_2^0,t_2)         = C(t_2)                                         & \Rightarrow & \breve C(t_2) = C(t_2) F(t_2,t_2^0).
\end{array}\right.
\]
Thus we obtain a differentiable vessel
\[ \begin{array}{ll}
\mathfrak{DV} = (F(t_2^0,t_2) A_1(t_2) F(t_2,t_2^0), I,F(t_2^0,t_2) \widetilde B(t_2),C(t_2) F(t_2,t_2^0), D(t_2),\widetilde D(t_2); \\
~~~~~~~~~~~~~~~~~~\sigma_1(t_2), \sigma_2(t_2), \gamma(t_2), \sigma_{1*}(t_2), \sigma_{2*}(t_2)\gamma_*(t_2);
        \mathcal H, \mathcal{E},\mathcal E_*, \mathcal{\widetilde E},\mathcal{\widetilde E_*})
\end{array} \]
with trivial evolution semi-group $I$, which is obviously absolutely continuous.
\qed

Notice that this vessel is of a very special form. Its evolution semi-group is trivial and as a result, differentiating this vessel,
we shall obtain that the generator of the semi group $A_2(t_2) = 0$ is trivial. 
Consequently, from the Lax equation, $A_1(t_2) = A_1$ becomes a constant operator.

Let us for the completeness of presentation explicitly write down the notion of quasi similarity of two differential vessels. Geven two
vessels
\[ \mathfrak{DV} = (A_1(t_2), A_2(t_2),\widetilde{B}(t_2),C(t_2), D(t_2),\widetilde D(t_2);
        \sigma_1(t_2), \sigma_2(t_2), \gamma(t_2), \sigma_{1*}(t_2), \sigma_{2*}(t_2),\gamma_*(t_2);
        \mathcal H, \mathcal{E},\mathcal E_*, \mathcal{\widetilde E},\mathcal{\widetilde E_*})
\]
and
\[ \breve{\mathfrak{DV}} = (\breve A_1(t_2), \breve A_2(t_2),\breve{\widetilde{B}}(t_2),\breve C(t_2), D(t_2),\widetilde D(t_2);
        \sigma_1(t_2), \sigma_2(t_2), \gamma(t_2), \sigma_{1*}(t_2), \sigma_{2*}(t_2),\gamma_*(t_2);
        \mathcal H, \mathcal{E},\mathcal E_*, \mathcal{\widetilde E},\mathcal{\widetilde E_*})
\]
we will say that they are quasi-similar if there exists a (possibly unbounded) linear operator $T(t_2): D(T(t_2)) \rightarrow \breve{\mathcal{H}}_{t_2}$, with a dense domain $D(T(t_2)) \subseteq \mathcal H_{t_2}$, which is 1-1, with dense range, absolutely
continuous (in the norm operator topology) and satisfies the following intertwining conditions
\begin{equation} \label{eq:DUconnect} \left\{\begin{array}{llll}
\breve A_1(t_2) T(t_2)       & = T(t_2) A_1(t_2), \\
\breve A_2(t_2) T(t_2)       & = T(t_2) A_2(t_2) + \frac{d}{dt_2} T(t_2), \\
\widetilde {\breve B}(t_2)   & = T (t_2) \widetilde B(t_2), \\
\breve C(t_2) T(t_2)         & = C(t_2).
\end{array}\right.
\end{equation}
We shall also demand that
\begin{equation} 
\widetilde{\mathcal C}_{t_2} \subseteq D(T(t_2)), ~~~~~~~~ \operatorname{Im}A_2(t_2) \subseteq D(T(t_2)).
\end{equation}
If $T(t_2)$ is an invertible bounded operator, we shall say that two differential vessels are \textit{similar}.
\subsection{\label{sec:SepVar}Separation of variables and the notion of transfer function}
One of the reasons why overdetermined systems invariant in one direction are interesting is the possiblity tp perform a partial 
separation of variables.
Taking all the trajectory data in the form
\[ \begin{array}{lll}
u(t_1,t_2) = u_\lambda(t_2) e^{\lambda t_1}, \\
x(t_1,t_2) = x_\lambda(t_2) e^{\lambda t_1}, \\
y(t_1,t_2) = y_\lambda(t_2) e^{\lambda t_1},
\end{array} \]
we arrive at the notion of a transfer function. Note that $u(t_1,t_2), y(t_1,t_2)$ satisfy PDEs, but
$u_\lambda(t_2), y_\lambda(t_2)$ are solutions of ODEs with a spectral parameter $\lambda$,
\begin{eqnarray}
\label{eq:InCC} \lambda \sigma_2(t_2)  u_\lambda(t_2) - \sigma_1(t_2) \frac{\partial}{\partial t_2}u_\lambda(t_2) +
\gamma(t_2) u_\lambda(t_2) = 0, \\
\label{eq:OutCC} \lambda \sigma_{2*}(t_2) y_\lambda(t_2) - \sigma_{1*}(t_2) \frac{\partial}{\partial t_2}y_\lambda(t_2) + \gamma_*(t_2)
y_\lambda(t_2) = 0.
\end{eqnarray}
The corresponding i/s/o system becomes
\[ \left\{ \begin{array}{lll}
    x_\lambda(t_2) = (\lambda I - A_1(t_2))^{-1} \widetilde B(t_2) \sigma_1(t_2) u_\lambda(t_2) \\
    x(t_1, t_2) = F (t_2,\tau_2) x(t_1, \tau_2) + \int\limits_{\tau_2}^{t_2} F(t_2, s) \widetilde B(s) \sigma_2(s) u(t_1, s)ds \\
    y_\lambda(t_2) = D(t_2)u_\lambda(t_2) +C(t_2) x_\lambda(t_2)
    \end{array} \right.
\]
The output $y_\lambda(t_2) = D(t_2) u_\lambda(t_2) + C(t_2) x_\lambda(t_2)$ may be found from the first i/s/o equation:
\[ y_\lambda(t_2) = S(\lambda, t_2) u_\lambda(t_2), \]
using the \textit{transfer function}
\begin{equation} \label{eq:S}
    S(\lambda, t_2) = D(t_2) + C(t_2) (\lambda I - A_1(t_2))^{-1} \widetilde B(t_2) \sigma_1(t_2).
\end{equation}
Here $\lambda$ is outside the spectrum of $A_1(t_2)$, which is independent of $t_2$ by (\ref{eq:LaxCond}). We
emphasize that $S(\lambda,t_2)$ is a function of $t_2$ for each $\lambda$ (which is a frequency
variable corresponding to $t_1$).
\begin{thm} \label{thm:PropS} A transfer function $S(\lambda,t_2)$, defined by (\ref{eq:S}) has the following properties:
    \begin{enumerate}
			\item For almost all $t_2$, $S(\lambda, t_2)$ is an analytic function of $\lambda$
				in the neighborhood of $\infty$, where it satisfies:
				\[ S(\infty, t_2) = D(t_2) \].
			\item For all $\lambda$, $S(\lambda, t_2)$ is an absolutely continuous function of $t_2$.
			\item For each fixed $\lambda$, multiplication by $S(\lambda, t_2)$ maps
				solutions of \linebreak $\lambda \sigma_2(t_2) u - \sigma_1(t_2) \frac{d u}{dt_2} + \gamma(t_2) u = 0$ to solutions of \newline
				$\lambda \sigma_2(t_2) y - \sigma_1(t_2) \frac{d y}{dt_2} + \gamma_*(t_2) y = 0$.
    \end{enumerate}
\end{thm}
\textbf{Proof:}
Those are easily checked properties, following from the definition of $S(\lambda,t_2)$:
\[ S(\lambda, t_2) = D(t_2) + C(t_2) (\lambda I - A_1(t_2))^{-1} \widetilde B(t_2) \sigma_1(t_2). \]
When $\lambda \rightarrow \infty$, since all the operators are
bounded the second summand vanishes and we obtain $ S(\infty, t_2)
= D(t_2)$. Moreover, it will be an analytic function of
$\lambda$, when $\lambda > \| A_1(t_2)\|$ and we obtain the first
property.

In order to understand the second property let us rewrite $S(\lambda, t_2)$, using the Lax equation in the following way:
\[ \begin{array}{lll}
    S(\lambda, t_2) & = D(t_2) + C(t_2) (\lambda I - A_1(t_2))^{-1} B(t_2) \sigma_1(t_2) = \\
       & = D(t_2) + C(t_2) (\lambda I - F(t_2,t_2^0)A_1(t_2^0) F(t_2^0, t_2))^{-1} B(t_2) \sigma_1(t_2) = \\
       & = D(t_2) + C(t_2) F(t_2,t_2^0)(\lambda I - A_1(t_2^0))^{-1} F(t_2^0, t_2) B(t_2) \sigma_1(t_2).
\end{array} \]
The functions $C(t_2) F(t_2,t_2^0)$, $F(t_2^0, t_2) B(t_2)$,
$\sigma_1(t_2)$ are absolutely continuous in appropriate spaces,
thus their multiplication too and we obtain the second property.

The third property is a direct result of our construction. \qed

\noindent\textbf{Remark:} for the case $\dim\mathcal H_{t_2}<\infty$, we obtain that $S(\lambda,t_2)$ is a rational (off the spectrum of $A_1(t_2)$)
in $\lambda$ function for all $t_2$.
\subsection{Class $\boldsymbol{\mathcal I}$ of intertwining functions}
We saw in previous section (thorem \ref{thm:PropS}) that transfer functions are very natural objects to study and have
three important properties. Suppose that we start from two ODEs the input \ref{eq:InCC} and the output \ref{eq:OutCC}
and denote the fundamental solution for them by $\Phi(\lambda,t_2, t_2^0)$ and by $\Phi_*(\lambda,t_2, t_2^0)$ respectively, then
\begin{equation} \label{eq:SInttw}
S(\lambda, t_2) \Phi(\lambda,t_2,t_2^0) = \Phi_*(\lambda,t_2,\tau_2) S(\lambda, t_2^0)
\end{equation}
and $S(\lambda,t_2)$ satisfies the following ODE
\begin{equation} \label{eq:DforS}
\frac{\partial}{\partial t_2} S(\lambda,t_2) = \sigma_{1*}^{-1}(t_2) (\sigma_{2*}(t_2) \lambda + \gamma_*(t_2)) S(\lambda,t_2)-
S(\lambda,t_2)\sigma_1^{-1}(t_2) (\sigma_2(t_2) \lambda + \gamma(t_2)).
\end{equation}
We recall \cite{bib:CoddLev} that from the fundamental theory of linear differential equations that for each equation there correspond and
invertible matrix (or operator) function $\phi(t_2,t_2^0)$ which obtains value $I$ for a fixed value of $t_2 = t_2^0$ and any other solution 
$u(t_2)$, satisfying $u(t_2^0) = u_0$ is just of the form 
\[ u(t_2) = \phi(t_2,t_2^0) u_0.
\]

Some of the simple properties of the fundamental matrix $\Phi(\lambda,t_2, t_2^0)$ are as follows. Notice that $\Phi(\lambda,t_2,t_2^0)$ can be replaced by $\Phi_*(\lambda,t_2,t_2^0)$ with a corresponding change of operators.
\begin{lemma} The fundamental matrix $\Phi(\lambda,t_2,t_2^0)$ in (\ref{eq:InCC}) satisfies:
\begin{eqnarray}
\frac{\partial}{\partial t_2} \Phi(\lambda,t_2,t_2^0) = \sigma_1^{-1}(\lambda \sigma_2(t_2) + \gamma(t_2)) \Phi(\lambda,t_2,t_2^0), \\
\label{eq:PhiInverse}
\frac{\partial}{\partial t_2} \Phi^{-1}(\lambda,t_2,t_2^0) = - \Phi^{-1}(\lambda,t_2,t_2^0) \sigma_1^{-1}(\lambda \sigma_2(t_2) + \gamma(t_2)), \\
\frac{\partial}{\partial t_2} \Phi^*(\lambda,t_2,t_2^0) = \Phi^*(\lambda,t_2,t_2^0) (\bar\lambda \sigma_2^*(t_2) + \gamma^*(t_2))\sigma_1^{-1*}, \\
\frac{\partial}{\partial t_2} \Phi^{-1*}(\lambda,t_2,t_2^0) =
        -(\bar\lambda \sigma_2^*(t_2) + \gamma^*(t_2)) \sigma_1^{-1*} \Phi(\lambda,t_2,t_2^0)^{-1*},
\end{eqnarray}
\end{lemma}
\textbf{Proof:} Conjugating (\ref{eq:InCC}) and using a formula for the derivative of the inverse. \qed

Since we will intensively work with such functions, we define their class $\boldsymbol{\mathcal I}$ as follows
\begin{defn} The class
\[ \boldsymbol{\mathcal I} = \boldsymbol{\mathcal I}(\sigma_1,\sigma_2, \gamma,\sigma_{1*},\sigma_{2*},\gamma_*)
\]
is a class of functions $S(\lambda,t_2)$ of two variables, which are
\begin{enumerate}
	\item analytic in the neighborhood of $\lambda=\infty$ for all $t_2$,
	\item absolutely continuous as functions of $t_2$ for almost all $\lambda$,
	\item map solutions of the input ODE (\ref{eq:InCC}) with spectral parameter $\lambda$ to the output 
		ODE (\ref{eq:OutCC}) with the same spectral parameter (i.e., they satisfy the equation (\ref{eq:SInttw}))
\end{enumerate}
\end{defn}

\section{\label{sec:MTGSim}Main theorem of gauge quasi similarity}
The following result is an analogue in our framework of the standard quasi similarity theorem for minimal systems
\cite{bib:Helton, bib:BallCohen}.
\begin{thm}\label{mthm:GE}
Assume that we are given two minimal (integral) vessels $\mathfrak{IV}, \breve{\mathfrak{IV}}$
\[ \begin{array}{ll}
\mathfrak{IV}=(A_1(t_2),F(t_2,t_2^0),\widetilde{B}(t_2),C(t_2),D(t_2),\widetilde
D(t_2);
        \sigma_1(t_2), \sigma_2(t_2), \gamma(t_2), \sigma_{1*}(t_2), \sigma_{2*}(t_2)\gamma_*(t_2);
        \mathcal H_{t_2}, \mathcal{E},\mathcal E_*, \mathcal{\widetilde E},\mathcal{\widetilde E_*}) \\
\mathfrak{\breve IV} = (\breve A_1(t_2), \breve
F(t_2,t_2^0),\widetilde{\breve B}(t_2),\breve C(t_2),
\breve D(t_2),\breve{\widetilde D}(t_2);
        \sigma_1(t_2), \sigma_2(t_2), \gamma(t_2), \sigma_{1*}(t_2), \sigma_{2*}(t_2)\gamma_*(t_2);
        \mathcal{\breve H}_{t_2}, \mathcal{E},\mathcal E_*, \mathcal{\widetilde E},\mathcal{\widetilde E_*})
\end{array} \]
with transfer functions $S(\lambda,t_2)$, $\breve{S}(\lambda,t_2)$. Then the vessels are gauge quasi similar iff
$S(\lambda,t_2) = \breve{S}(\lambda,t_2) $ in a neighborhood of $\lambda=\infty$.
\end{thm}
\textbf{Proof:} The easy direction of the statement considers the case when there exists $T(t_2):\mathcal H_{t_2} \rightarrow 
\breve{\mathcal H}_{t_2}$, responsible for quasi-similarity, i.e., satisfying (\ref{eq:Uconnect}) and (\ref{eq:TDomainCond}),
then
\[ \begin{array}{llllll}
S(\lambda,t_2) & = D(t_2) +  C(t_2)(\lambda I - A_1(t_2))^{-1} \widetilde{B}(t_2) \sigma_1(t_2) =
D(t_2) +  \breve C(t_2) T(t_2) (\lambda I - A_1(t_2))^{-1} \widetilde{B}(t_2) \sigma_1(t_2) = \\
& = D(t_2) +  \breve C(t_2) (\lambda I - \breve A_1(t_2))^{-1} T(t_2) \widetilde{B}(t_2) \sigma_1(t_2) =
D(t_2) +  \breve C(t_2) (\lambda I - \breve A_1(t_2))^{-1} \widetilde{\breve B}(t_2) \sigma_1(t_2) = \\
& = \breve S(\lambda,t_2).
\end{array} \]
Fo the converse, we obtain first that the values at infinity of the two functions are equal:
\[ D(t_2) = \breve D(t_2) \Rightarrow  \widetilde D(t_2) = \breve{\widetilde D}(t_2) \text{, since $\sigma_1, \sigma_{1*}$ are invertible.}
\]
By looking at the Taylor coefficients in the power series expansions of functions $S(\lambda,t_2)$,
$\breve {S}(\lambda,t_2)$ at infinity, we obtain
\begin{equation} \label{eq:ExpEqual}
 C(t_2) A_1^n(t_2) \widetilde B(t_2) = \breve C(t_2) \breve A_1^n(t_2) \widetilde {\breve B}(t_2), \text{ for all $n = 0,1,2,\ldots$}
\end{equation}
Since $\mathfrak{V}$ is approximately controllable, the set
\[ \widetilde{\mathcal C}_{t_2} = \{ \sum_{i=1}^N F(t_2,s_i) A_1^{n_i}(s_i) \widetilde B(s_i) e_i \mid
        n_i, N \in \mathbb N, e_i\in \mathcal E, s_i \in \mathbb R\}
\]
is dense in $\mathcal H_{t_2}$. Define next $T: \widetilde{\mathcal C}_{t_2} \rightarrow \mathcal{\breve H}_{t_2}$ by
\[ T (\sum_{i=1}^N F(t_2,s_i) A_1^{n_i}(s_i) \widetilde B(s_i) e_i) =
                \sum_{i=1}^N \breve F(t_2,s_i) \breve A_1^{n_i}(s_i) \widetilde{\breve B}(s_i) e_i.
\]
Then $T$ is obviously a linear transformation, provided it is well defined, i.e. we have to check that
\[ 0 = \sum_{i=1}^N F(t_2,s_i) A_1^{n_i}(s_i) \widetilde B(s_i) e_i \Rightarrow 0 =
                \sum_{i=1}^N \breve F(t_2,s_i) \breve A_1^{n_i}(s_i) \widetilde{\breve B}(s_i) e_i.
\]
Since by assumption $\breve{\mathfrak{V}}$ is observable, to show
$0 = \sum_{i=1}^N \breve F(t_2,s_i) \breve A_1^{n_i}(s_i) \widetilde{\breve B}(s_i) e_i$ is the same as showing that
\[ \breve C(t_2) \breve A_1^k(t_2) (\sum_{i=1}^N \breve F(t_2,s_i) \breve A_1^{n_i}(s_i) \widetilde{\breve B}(s_i) e_i)  = 0 , \text{ for all $k = 0,1,2,\ldots$}
\]
This is done with the help of the following lemma
\begin{lemma} For each $t_2, s$ the following equality holds:
\[ C(t_2) (\lambda I -A_1(t_2))^{-1} F(t_2,s) \widetilde B(s) \sigma_1(s) e =
\breve C(t_2) (\lambda I - \breve A_1(t_2))^{-1} \breve F(t_2,s) \widetilde{\breve B}(s) \sigma_1(s) e
\]
\end{lemma}
\textbf{Proof:} Simple calculations using the first vessel condition \ref{eq:OverDetCondIn}, result in
\[ \begin{array}{llrr}
\frac{\partial}{\partial s} \big( C(t_2) (\lambda I -A_1(t_2))^{-1} F(t_2,s) \widetilde B(s) \sigma_1(s) \big) = \\
= - C(t_2) (\lambda I -A_1(t_2))^{-1} F(t_2,s) \widetilde B(s) [\sigma_2(s) \lambda + \gamma(s)] + C(t_2) F(t_2,s) \widetilde B(s) \sigma_2(s),
\end{array} \]
Remember that the function $\Phi(\lambda,t_2,t_2^0)$ satisfies (\ref{eq:PhiInverse}), so using variation of coefficients
\[ C(t_2) (\lambda I -A_1(t_2))^{-1} F(t_2,s) \widetilde B(s) \sigma_1(s) = K(\lambda,t_2,s) \Phi^{-1} (\lambda,t_2,s)
\]
where
\[ \frac{\partial}{\partial s} K(\lambda,t_2,s) = C(t_2) F(t_2,s) \widetilde B(s) \sigma_2(s) \Phi (\lambda,t_2,s).
\]
Particularly, $K(\lambda,t_2,s)$ is 
\begin{multline*}
K(\lambda,t_2,s) = K(\lambda,t_2^0,s) +\int_{t_2^0}^{t_2}  \frac{\partial}{\partial s} K(\lambda,t_2,s) ds = \\
= K(\lambda,t_2^0,s) +\int_{t_2^0}^{t_2}  C(t_2) F(t_2,s) \widetilde B(s) \sigma_2(s) \Phi (\lambda,t_2,s) ds
\end{multline*}
and all its poles are at the first term $K(\lambda,t_2^0,s)$. The same considerations, applied to
$ \breve C(t_2) (\lambda I - \breve A_1(t_2))^{-1} \breve F(t_2,s) \widetilde{\breve B}(s) \sigma_1(s)$
results in a function $\breve K(\lambda,t_2,s)$, whose poles are also at the term $\breve K(\lambda,t_2^0,s)$:
\[ \breve K(\lambda,t_2,s) = \breve K(\lambda,t_2^0,s) +\int_{t_2^0}^{t_2} \breve C(t_2) \breve F(t_2,s) \widetilde{\breve B}(s) \sigma_2(s)
\Phi (\lambda,t_2,s) ds. \]
Then
\[ \begin{array}{llrr}
C(t_2) (\lambda I -A_1(t_2))^{-1} F(t_2,s) \widetilde B(s) \sigma_1(s) -
        \breve C(t_2) (\lambda I - \breve A_1(t_2))^{-1} \breve F(t_2,s) \widetilde{\breve B}(s) \sigma_1(s)  = \\
= [K(\lambda,t_2,s) - \breve K(\lambda,t_2,s)] \Phi^{-1} (\lambda,t_2,s) = \\
= [K(\lambda,t_2,t_2) + \int_{t_2}^s \frac{\partial}{\partial y}K(\lambda,t_2,y) dy- \breve K(\lambda,t_2,t_2) -
        \int_{t_2}^s \frac{\partial}{\partial y} \breve K(\lambda,t_2,y) dy] \Phi^{-1} (\lambda,t_2,s) = \\
= [\int_{t_2}^s \frac{\partial}{\partial y}K(\lambda,t_2,y) dy - \int_{t_2}^s \frac{\partial}{\partial y}\breve K(\lambda,t_2,y) dy]\Phi^{-1} (\lambda,t_2,s),
\end{array} \]
in other words, this difference is an entire function of $\lambda$. On the other hand, two functions
\[ C(t_2) (\lambda I -A_1(t_2))^{-1} F(t_2,s) \widetilde B(s) \sigma_1(s), ~~~
\breve C(t_2) (\lambda I - \breve A_1(t_2))^{-1} \breve F(t_2,s) \widetilde{\breve B}(s) \sigma_1(s)
\]
are zero at infinity (i.e., globally bounded). Taking a paring with an arbitrary linear functional
and the operator applied to arbitrary vector, by Liouville’s theorem this difference is constant and is equal to the value at infinity for
each such pairing. Thus the operator itself is zero. \qed

One of the consequences of this theorem is that Taylor coefficients around infinity of the functions
\[ C(t_2) (\lambda I -A_1(t_2))^{-1} F(t_2,s) \widetilde B(s) \sigma_1(s), ~~~
\breve C(t_2) (\lambda I - \breve A_1(t_2))^{-1} \breve F(t_2,s) \widetilde{\breve B}(s) \sigma_1(s)
\]
are equal, and consequently,
\[ C(t_2) A_1(t_2)^n F(t_2,s) \widetilde B(s) \sigma_1(s) =
        \breve C(t_2) \breve A_1(t_2)^n \breve F(t_2,s) \widetilde{\breve B}(s) \sigma_1(s)
\]
Using this result, the fact that $T(t_2)$ is well defined is immediate, because
\[ \begin{array}{llll}
\breve C(t_2) \breve A_1^k(t_2) (\sum_{i=1}^N \breve F(t_2,s_i) \breve A_1^{n_i}(s_i) \widetilde{\breve B}(s_i) e_i) =
\sum\limits_{i=0}^N \breve C(t_2) \breve F(t_2,s_i) \breve A_1^{n_i+k}(s_i) \widetilde{\breve B}(s_i) e_i = \\
~~ = \sum\limits_{i=0}^N C(t_2) F(t_2,s_i) A_1^{n_i+k}(s_i) \widetilde{B}(s_i) e_i =
C(t_2) A_1^k(t_2) (\sum\limits_{i=0}^N F(t_2,s_i) A_1^{n_i}(s_i) \widetilde{B}(s_i) e_i = C(t_2) A^k(t_2) 0 = 0.
\end{array} \]
One checks that $T$ is one-to-one by using the observability of $\mathfrak{V}$:
\[ \begin{array}{llll}
\operatorname{Ker} T = \bigcap_{n\in\mathbb N, s\in\mathbb R} \operatorname{Ker} C(s) A_1^n(s) F(s,t_2) = \widetilde{\mathcal O}^\perp_{t_2},
\end{array} \]
which is trivial by assumption of minimality. That $T$ has dense range follows from the approximate controllability of $\breve{\mathfrak{V}}$. 
This finishes the proof. \qed

Since the notions of similarity and quasi-similarity are identical for finite dimensional vessels, we obtain
\begin{cor} Two finite dimensional vessels (with rational in $\lambda$ transfer functions) are gauge similar,
if and only if, the transfer functions are identical at the neighborhood of $\lambda=\infty$.
\end{cor}
\section{\label{sec:Adj}Adjoint system}
The notion of the adjoint system is very useful in system theory. It is obtained from a simple observation
that applying adjoint to the vessel conditions gives rise to a new set of conditions on
adjoint operators, which are almost vessel conditions. Moreover $(\Sigma^*)^*$ is actually
$\Sigma$ by a trivial change of coordinates $x \rightarrow -x$ on the state space. Here is the precise definition.
Given a system in the differential form
\begin{equation} \label{eq:DSigma}
    \mathfrak D\Sigma: \left\{ \begin{array}{lll}
    \frac{\partial}{\partial t_1}x(t_1,t_2) = A_1(t_2) x(t_1,t_2) + \widetilde{B}(t_2) ~\sigma_1(t_2) ~u(t_1,t_2) \\
    \frac{\partial}{\partial t_2}x(t_1,t_2) = A_2(t_2) x(t_1,t_2) + \widetilde{B}(t_2) ~\sigma_2(t_2) ~u(t_1,t_2) \\
    y(t_1,t_2) = D(t_2) u(t_1,t_2) + C(t_2) x(t_1,t_2)
    \end{array} \right.
\end{equation}
and associated vessel
\[ \mathfrak{V} = (A_1(t_2), A_2(t_2),\widetilde{B}(t_2),C(t_2), D(t_2),\widetilde D(t_2);
        \sigma_1(t_2), \sigma_2(t_2), \gamma(t_2), \sigma_{1*}(t_2), \sigma_{2*}(t_2)\gamma_*(t_2);
        \mathcal H, \mathcal{E},\mathcal E_*, \mathcal{\widetilde E},\mathcal{\widetilde E_*})
\]
it is natural to introduce the \textit{adjoint system}:
\begin{equation} \label{eq:systemAdj}
    D\Sigma^*: \left\{ \begin{array}{lll}
    -\frac{\partial}{\partial t_1}x_*(t_1,t_2) = A_1^*(t_2) ~x_*(t_1,t_2) + C^*(t_2) \sigma^*_{1*} u_*(t_1,t_2) \\[5pt]
    -\frac{\partial}{\partial t_2}x_*(t_1,t_2) = A_2^*(t_2) ~x_*(t_1,t_2) + C^*(t_2) \sigma^*_{2*} u_*(t_1,t_2) \\[5pt]
    y_*(t_1,t_2) = \widetilde B^*~ x_*(t_1,t_2)  + \widetilde D^*(t_2)~ u_*(t_1,t_2),
    \end{array} \right.
\end{equation}
which is associated to the vessel $\mathfrak{V^*}$ given by
\[ \mathfrak{V^*} =
(-A_1^*, -A_2^*, -C^*, \widetilde B^*, \widetilde D^*, D^*;
\sigma^*_{1*}, \sigma^*_{2*}, -\gamma_*^* - \frac{d}{dt_2}
\sigma_{1*}^*, \sigma^*_1, \sigma^*_2, -\gamma^* - \frac{d}{dt_2}
\sigma_1^*; \mathcal{H}, \widetilde{\mathcal E}_*,
\widetilde{\mathcal E}_*, \mathcal{E}_*, \mathcal{E}),
\]
where all the operators are functions of $t_2$ and satisfy the following axioms:
\[ \begin{array}{lllllll}
    \frac{d}{dt_2} A_1^* = A_1^* A_2^* - A_2^* A_1^* \\
   -\frac{d}{dt_2} \big(C^* \sigma^*_{1*} \big) - A_2^* C^* \sigma^*_{1*} + A_1^* C^* \sigma^*_{2*}
   + C^* (\gamma_*^*+\frac{d}{dt_2} \sigma_{1*}^*) = 0 \\
    \frac{d}{dt_2} \big(\sigma_1^* \widetilde B^* \big) - \sigma_1^* \widetilde B^* A_2^* +
    \sigma_2^* \widetilde B A_1^* + \gamma^* \widetilde B^* = 0 \\
   \sigma^*_1 \widetilde D^* = D^* \sigma^*_{1*}, ~~ \sigma^*_2 \widetilde D^* = D^* \sigma^*_{2*} \\
   D^* (-\gamma^*_* -\frac{d}{dt_2} \sigma_{1*}^*) = -\sigma^*_2 \widetilde B^* C^* \sigma^*_{1*} +
   \sigma^*_1 \widetilde B^* C^* \sigma^*_{2*} - \sigma^*_1 \frac{d}{dt_2} \widetilde D^* - (\gamma^* +\frac{d}{dt_2} \sigma_1^*)\widetilde D^*.
   \end{array}
\]
Moreover, the transfer function of the adjoint vessel $S_*(\mu,t_2)$ maps solutions of the \textit{adjoint input} ODE
\begin{equation} \label{eq:InCC*} [\sigma^*_{2*} \mu - \sigma^*_{1*} \dfrac{d}{dt_2} -
\gamma^*_* - \frac{d}{dt_2} \sigma_{1*}^*] u_*(t_2) = 0
\end{equation}
with the spectral parameter $\mu$ to solutions of the \textit{adjoint output} ODE
\begin{equation} \label{eq:OutCC*} [\sigma^*_2 \mu - \sigma^*_1 \dfrac{d}{dt_2} -
\gamma^*-\frac{d}{dt_2} \sigma_1^*] y_*(t_2) = 0
\end{equation}
with the same spectral parameter. In the language of fundamental
matrices, if one denotes the fundamental solution for $u_*(t_2),
y_*(t_2)$ as $\Psi(\mu,t_2, t_2^0)$ and $\Psi_*(\mu,t_2, t_2^0)$,
respectively, then (similarly to (\ref{eq:SInttw}))
\[ S_*(\mu, t_2) \Psi(\mu,t_2,\tau_2) =
\Psi_*(\mu,t_2,\tau_2) S_*(\mu, \tau_2)
\]
Notice also that the equation for $S_*(\lambda, t_2)$ similar to
(\ref{eq:S}) is:
\[ S_*(\mu, t_2) = \widetilde D^* - \widetilde B^* (\mu I +
A_1^*)^{-1} C^* \sigma^*_{1*}
\]
And as in the case of constant operators, one obtains that
$\mathfrak{V^*}$ is a vessel iff $\mathfrak{V}$ is. And as we mentioned before,
$\mathfrak{V^*}^*$ is the same as $\mathfrak V$ after a trivial
change of coordinates $x \rightarrow -x$ on the state space.

\noindent \textbf{Remarks:}
\begin{enumerate}
\item It is easy to check that $\sigma_{1*}^* \Phi_*^{-1*}(\lambda, t_2, t_2^0)$ and
$\Psi(-\bar\lambda, t_2, t_2^0)$ satisfy the same differential equation. Since they are fundamental
matrices it is possible iff
\begin{eqnarray}\label{eq:Phi*Psi}
 \sigma_{1*}^{-1*} \Phi_*^{-1*}(\lambda, t_2, t_2^0) = \Psi(-\bar\lambda, t_2, t_2^0)  \sigma_{1*}^{-1*} \\
 \label{eq:PhiPsi*} \sigma_1^{-1*} \Phi^{-1*}(\lambda, t_2, t_2^0) = \Psi_*(-\bar\lambda, t_2, t_2^0)  \sigma_1^{-1*}
\end{eqnarray}
\item This relation between $\Phi$ and $\Psi$ fundamental matrices means that the following relation
between transfer functions has to be satisfied:
\begin{equation} \label{eq:SS*}
 S(\lambda, t_2) =  \sigma_{1*}^{-1} S^*_*(-\bar\lambda, t_2) \sigma_1
\end{equation}
One can easily verify this formula directly, using the vessel conditions and the formulas for
$S(\lambda, t_2)$, $S^*_*(-\bar\lambda, t_2)$:
\[ \begin{array}{lllllll}
  \sigma_{1*}^{-1} S^*_*(-\bar\lambda, t_2) \sigma_1 =
  \sigma_{1*}^{-1} \big[ \widetilde D^* - \widetilde B^* (-\bar\lambda I +
A_1^*)^{-1} C^* \sigma^*_{1*} \big]^* \sigma_1 = \\
 ~~~~ = \sigma_{1*}^{-1} [\widetilde D - \sigma_{1*} C(-\lambda I + A_1)^{-1} \widetilde B]\sigma_1 =
    \sigma_{1*}^{-1} \widetilde D \sigma_1 - C(-\lambda I + A_1)^{-1} \widetilde B \sigma_1 = \\
 ~~~~ = D + C (\lambda I - A_1)^{-1} \widetilde B \sigma_1 = S(\lambda, t_2).
\end{array} \]
\item For values of $\lambda$ outside of the spectrum of $A_1(t_2)$ (which is independent of $t_2$),
we claim that $S(\lambda, t_2)$ is invertible iff $S^*_*(-\bar\lambda, t_2)$ is and from (\ref{eq:SS*})
we conclude that
\[ S(\lambda, t_2)^{-1} = \sigma_1^{-1} S^{-1*}_*(-\bar\lambda, t_2) \sigma_{1*} \]
\end{enumerate}

\section{\label{sec:Operations}System and vessel operations}
As in the classical case we develop basic system operations in our setting, which will be of great importance in solving
factorization problems for transfer functions.
\subsection{Cascade connection of systems}
Suppose we are given two vessels
\begin{equation}\label{eq:TwoCoupledVessels} \begin{array}{llll}
\mathfrak{IV}^\prime = (A_1^\prime(t_2), F^\prime(t_2,t_2^0),\widetilde B^\prime(t_2), C^\prime(t_2), D^\prime(t_2),\widetilde D^\prime(t_2); \\
~~~~~~~~~~~~~~~~~~~~~~~~~~~~~~~~~~~~~~~~~~~~~~~~\sigma'_1(t_2), \sigma'_2(t_2), \gamma'(t_2), \sigma'_{1*}(t_2), \sigma'_{2*}(t_2)\gamma'_*(t_2); \mathcal H^\prime_{t_2}, \mathcal{E}',\mathcal E'_*, \mathcal{\widetilde E'},\mathcal{\widetilde E'_*}), \\
\mathfrak{IV}^{\prime\prime} = (A_1^{\prime\prime}(t_2), F^{\prime\prime}(t_2,t_2^0),\widetilde B^{\prime\prime}(t_2), C^{\prime\prime}(t_2), D^{\prime\prime}(t_2),\widetilde D^{\prime\prime}(t_2); \\
~~~~~~~~~~~~~~~~~~~~~~~~~~~~~~~~~~~~~~~~~~~~~~~~\sigma''_1(t_2), \sigma''_2(t_2), \gamma''(t_2),
\sigma''_{1*}(t_2), \sigma''_{2*}(t_2)\gamma''_*(t_2); \mathcal H^{\prime\prime}_{t_2}, \mathcal{E}'',\mathcal E''_*,\mathcal{\widetilde E''},\mathcal{\widetilde E''_*}).
\end{array} \end{equation}
and the corresponding systems $I\Sigma^{\prime}, I\Sigma^{\prime\prime}$, defined in (\ref{eq:system}).
We want to generate a new system $I\Sigma$ by feeding in the output of the first system $I\Sigma^{\prime}$ as the input for
the second system $I\Sigma^{\prime\prime}$. To this end we assume the output spaces of the first system are the same as the input spaces of
the second system (as in the case of the classical cascade connection) but also that {\em the corresponding compatibility conditions
hold}:
\begin{equation} \label{eq:CascConComp}
 \sigma'_{1*}(t_2) = \sigma''_1(t_2), ~~\sigma'_{2*}(t_2)= \sigma''_2(t_2), ~~\gamma'_*(t_2) = \gamma''(t_2),~~~~~~~
\mathcal{\widetilde E'} = \mathcal{E}'', ~~\mathcal{\widetilde E'_*} = \mathcal E''_*
\end{equation}
Thus we obtain the following system of equations
\[ \begin{array}{lll}
I\Sigma': \left\{ \begin{array}{lll}
    \frac{\partial}{\partial t_1}x'(t_1,t_2) = A_1^\prime(t_2) x'(t_1,t_2) + \widetilde B^\prime(t_2) \sigma'_1(t_2) u(t_1,t_2) \\
    x'(t_1, t_2) = F^\prime (t_2,t_2^0) x'(t_1, t_2^0) + \int\limits_{t_2^0}^{t_2} F^\prime(t_2, s) \widetilde B^\prime(s) \sigma'_2(s) u(t_1, s)ds \\
    y'(t_1, t_2) = C^\prime(t_2) x'(t_1, t_2) + D^\prime(t_2) u(t_1, t_2)
  \end{array} \right. \\
~~~~~~~~~~~~~~~~~~~~~~~~~~~~~~~~~~~~~ I\Sigma'': \left\{ \begin{array}{lll}
    \frac{\partial}{\partial t_1}x''(t_1,t_2) = A''_1(t_2) x''(t_1,t_2) + \widetilde B''(t_2)\sigma'_{1*}(t_2) u''(t_1,t_2) \\
    x''(t_1, t_2) = F''(t_2,t_2^0) x''(t_1, t_2^0) + \int\limits_{t_2^0}^{t_2} F''(t_2, s) \widetilde B''(s)\sigma'_{2*}(s) u''(t_1, s)ds \\
    y(t_1, t_2) = C''(t_2) x''(t_1, t_2) + D''(t_2) u''(t_1, t_2)
  \end{array} \right.
\end{array} \]
Setting $u''(t_2) = y'(t_2)$, eliminating it and simplifying we get
\[
I\Sigma: \left\{ \begin{array}{lllll}
    \frac{\partial}{\partial t_1}\bbmatrix{x'(t_1,t_2)\\x''(t_1,t_2)}  = \bbmatrix{A'_1(t_2)& 0 \\\widetilde B''(t_2) \sigma'_1(t_2)C'(t_2) & A_1''(t_2)}
    \bbmatrix{x'(t_1,t_2)\\x''(t_1,t_2)} + \bbmatrix{B'(t_2)\\ \widetilde B''(t_2) \widetilde D'(t_2)} \sigma'_1 u(t_1,t_2) \\
    \bbmatrix{x'(t_1,t_2)\\x''(t_1,t_2)} = \bbmatrix{F' (t_2,t_2^0) & 0 \\ \int_{t_2^0}^{t_2} F''(t_2,s) \widetilde B''(s)\sigma'_2(s) C'(s) F('s,t_2^0) ds & F''(t_2,t_2^0)}
    \bbmatrix{x'(t_1,t_2^0)\\x''(t_1,t_2^0)} + \\
   ~~~~~~~~ + \int\limits_{t_2^0}^{t_2} \bbmatrix{F' (t_2,s) & 0 \\ \int_s^{t_2} F''(t_2,y) \widetilde B''(y)\sigma'_2(y) C'(y) F'(y,t_2^0) dy & F''(t_2,s)}
            \bbmatrix{\widetilde B'(s) \\ \widetilde B''(s)\widetilde D'(s)} \sigma'_2(s) u(t_1, s)ds \\
    y(t_1, t_2) = \bbmatrix{D''(t_2) C'(t_2) & C''(t_2)} \bbmatrix{x'(t_1,t_2^0)\\x''(t_1,t_2^0)} + D''(t_2) D'(t_2) u(t_1, t_2)
  \end{array} \right.
\]
Thus the corresponding vessel of this system is the following
\begin{equation}\label{eq:CascVessel} \begin{array}{llll}
\mathfrak{IV} = (\bbmatrix{A'_1(t_2)& 0 \\ \widetilde B''(t_2) \sigma'_1(t_2)C'(t_2) & A_1''(t_2)},
            \bbmatrix{F' (t_2,t_2^0) & 0 \\ \int_{t_2^0}^{t_2} F''(t_2,s) \widetilde B''(s)\sigma'_2(s) C'(s) F('s,t_2^0) ds & F''(t_2,t_2^0)}, \\
~~~~~~~~~~~~~~~~~~\bbmatrix{B'(t_2)\\ \widetilde B''(t_2) \widetilde D'(t_2)} ,\bbmatrix{D''(t_2) C'(t_2) & C''(t_2)} , D''(t_2) D'(t_2),
                \widetilde D''(t_2)\widetilde D'(t_2); \\
~~~~~~~~~~~~~~~~~~~~~~~~~~\sigma'_1(t_2), \sigma'_2(t_2), \gamma'(t_2), \sigma''_{1*}(t_2), \sigma''_{2*}(t_2), \gamma''_*(t_2);
            \mathcal H'_{t_2}\oplus\mathcal H''_{t_2}, \mathcal{E}',\mathcal E'_*, \mathcal{\widetilde E''},\mathcal{\widetilde E''_*})
\end{array}
\end{equation}
Last evaluation suggests that the system $I\Sigma$ with compatibility conditions $u'(t_1,t_2) = y''(t_1,t_2)$ indeed corresponds to a vessel
\begin{thm} \label{thm:CascConInt}
Given two vessels $\mathfrak{IV}', \mathfrak{IV}''$, defined in (\ref{eq:TwoCoupledVessels}) and satisfying compatibility conditions
(\ref{eq:CascConComp})
\[ \sigma'_{1*}(t_2) = \sigma''_1(t_2), ~~\sigma'_{2*}(t_2)= \sigma''_2(t_2), ~~\gamma'_*(t_2) = \gamma''(t_2),~~~~~~~
\mathcal{\widetilde E'} = \mathcal{E}'', ~~\mathcal{\widetilde E'_*} = \mathcal E''_*
\]
their cascade connection $\mathfrak{IV}$, defined in (\ref{eq:CascVessel}) is a vessel. Transfer functions
$S'(\lambda,t_2), S''(\lambda,t_2), S(\lambda,t_2)$ of the vessels $\mathfrak{IV}', \mathfrak{IV}'', \mathfrak{IV}$ respectivley
satisfy the following relation
\[ S(\lambda,t_2) =  S''(\lambda,t_2) S'(\lambda,t_2)
\]
\end{thm}
\textbf{Proof:} We have already seen that $\mathfrak{IV}$ is a vessel, provided $\mathfrak{IV}', \mathfrak{IV}''$ are. In order to see
the formula for transfer functions, we feed in the output $y'_\lambda(t_2)$ of the first system as the input $u''_\lambda(t_2)$ for
the second system (recall definitions in section \ref{sec:SepVar}). Then $y'_\lambda(t_2) = S(\lambda,t_2) u'_\lambda(t_2)$ is the input ($u''_\lambda(t_2) = y'_\lambda(t_2)$) for the second system and
\[ y''_\lambda(t_2) = S''(\lambda,t_2) u''_\lambda(t_2) = S''(\lambda,t_2) y'_\lambda(t_2) = S''(\lambda,t_2) S'(\lambda,t_2) u'_\lambda(t_2).
\]
We conclude that the transfer function $S(\lambda,t_2)$ for the composite system $I\Sigma$ is simply the product of the transfer functions of
the component systems:
\begin{equation} \label{eq:CascS} S(\lambda,t_2) = S'(\lambda,t_2) S(\lambda,t_2).
\end{equation}
\qed

The following theorem is the analogue of theorem \ref{thm:CascConInt} for differential vessels. It can be either deduced from theorem
\ref{thm:CascConInt} by differentiation, or the formulas can be obtained directly by writing the systems equations in the differential form
and the result established directly by algebraic manipulations
\begin{thm} Suppose that we are given two differential vessels
\[ \begin{array}{llll}
 \mathfrak{DV}' = (A_1'(t_2), A_2'(t_2),\widetilde{B}'(t_2),C'(t_2), D'(t_2),\widetilde D'(t_2);\\
~~~~~~~~~~~~~~~~~~~~~~~~~~~ \sigma'_1(t_2), \sigma'_2(t_2), \gamma'(t_2), \sigma'_{1*}(t_2), \sigma'_{2*}(t_2),\gamma'_*(t_2);
        \mathcal H', \mathcal{E}',\mathcal E'_*, \mathcal{\widetilde E'},\mathcal{\widetilde E'_*})
\end{array} \]
and
\[ \begin{array}{llll}
 \mathfrak{DV}'' = (A_1''(t_2), A_2''(t_2),\widetilde{B}''(t_2), C''(t_2), D''(t_2),\widetilde D''(t_2);\\
~~~~~~~~~~~~~~~~~~~~~~~~~~~ \sigma''_1(t_2), \sigma''_2(t_2),\gamma''(t_2), \sigma''_{1*}(t_2), \sigma''_{2*}(t_2), \gamma''(t_2), \gamma''_*(t_2);
        \mathcal H'', \mathcal{E''},\mathcal E''_*, \mathcal{\widetilde E}'',\mathcal{\widetilde E_*}'').
\end{array} \]
satisfying the compatibility conditions (\ref{eq:CascConComp}), then the following collection
\[ \begin{array}{llll}
 \mathfrak{DV} = (\bbmatrix{A'_1(t_2)& 0 \\ \widetilde B_1''(t_2)\sigma_1(t_2) C'(t_2) & A''_1(t_2)},
        \bbmatrix{A'_2(t_2)& 0 \\ \widetilde B_1''(t_2)\sigma_2(t_2)C'(t_2) & A''_2(t_2)}, \\
~~~~~~~~~~~~~~~~~~~\bbmatrix{\widetilde B'(t_2) \\ \widetilde B''(t_2) \widetilde D'(t_2)} ,\bbmatrix{D''(t_2) C'(t_2) & C''(t_2)} , D''(t_2) D'(t_2),
                \widetilde D''(t_2)\widetilde D'(t_2);\\
~~~~~~~~~~~~~~~~~~~~~~~~~~~~~~~~~~~~~~ \sigma'_1(t_2), \sigma'_2(t_2), \gamma'(t_2), \sigma''_{1*}(t_2), \sigma''_{2*}(t_2),
                         \gamma''_*(t_2)(t_2); \mathcal H'\oplus\mathcal H'', \mathcal{E}',\mathcal E'_*, , \mathcal{\widetilde E}'',\mathcal{\widetilde E_*''})
\end{array} \]
is a vessel called cascade connection of the vessels $\mathfrak{DV}', \mathfrak{DV}''$. The transfer functions of the
corresponding systems satisfy the formula (\ref{eq:CascS}).
\end{thm}

\subsection{Inversion of systems}
For the classical case \cite{bib:BGK}, if the feed through operator $D$ is invertible,
then one can define an inverse system having a transfer function equal to the reciprocal
of the transfer function of the original system. The analogue for $t_1$ invariant overdetermined $2D$ system
is as follows. Suppose that for the vessel in the differential form
\[ \mathfrak{DV} = (A_1(t_2),  A_2(t_2),\widetilde{B}(t_2), C(t_2), D(t_2),\widetilde D(t_2);
        \sigma_1(t_2), \sigma_2(t_2), \gamma(t_2), \sigma_{1*}(t_2), \sigma_{2*}(t_2)\gamma_*(t_2);
        \mathcal H, \mathcal{E},\mathcal E_*, \mathcal{\widetilde E},\mathcal{\widetilde E_*})
\]
both $D: \mathcal E \rightarrow \mathcal E_*$ and $\widetilde D:\widetilde{\mathcal E} \rightarrow \widetilde{\mathcal E_*}$ are invertible.
Then we may solve $u$ in terms of $y$ from the last system equation
\[ D\Sigma: \left\{ \begin{array}{lll}
    \frac{\partial}{\partial t_1}x(t_1,t_2) = A_1(t_2) x(t_1,t_2) + \widetilde{B}(t_2) ~\sigma_1(t_2) ~u(t_1,t_2) \\
    \frac{\partial}{\partial t_2}x(t_1,t_2) = A_2(t_2) x(t_1,t_2) + \widetilde{B}(t_2) ~\sigma_2(t_2) ~u(t_1,t_2) \\
    y(t_1,t_2) = D(t_2) u(t_1,t_2) + C(t_2) x(t_1,t_2)
    \end{array} \right.
\]
by
\[ u(t_1,t_2) = - D^{-1}(t_2) C(t_2) x(t_1,t_2) + D^{-1}(t_2) y(t_1,t_2)
\]
and plug it back to get a system $D\Sigma^\times$ having the property that $(u,x,y)$ is a trajectory for $D\Sigma$ if and only if
$(y,x,u)$ is a trajectory for $D\Sigma^\times$:
\[ D\Sigma^\times: \left\{ \begin{array}{lll}
    \frac{\partial}{\partial t_1}x(t_1,t_2) = (A_1(t_2)-\widetilde B\sigma_1 D^{-1}(t_2)C(t_2) )x(t_1,t_2) + \widetilde{B}(t_2) ~\sigma_1(t_2)D^{-1}(t_2) ~y(t_1,t_2) \\
    \frac{\partial}{\partial t_2}x(t_1,t_2) = (A_2(t_2)-\widetilde B\sigma_2 D^{-1}(t_2)C(t_2) ) x(t_1,t_2) + \widetilde{B}(t_2) ~\sigma_2(t_2)D^{-1}(t_2) ~y(t_1,t_2) \\
    u(t_1,t_2) = - D^{-1}(t_2) C(t_2) x(t_1,t_2) + D^{-1}(t_2) y(t_1,t_2)
    \end{array} \right.
\]
The linkage conditions (\ref{eq:LinkCond}) means that
\[ \widetilde D^{-1} \sigma_{1*} = \sigma_1 D^{-1}, ~~~~~~ \widetilde D^{-1} \sigma_{2*} = \sigma_2 D^{-1},
\]
so that this system can be rearranged somewhat to
\[ D\Sigma^\times: \left\{ \begin{array}{lll}
    \frac{\partial}{\partial t_1}x(t_1,t_2) = (A_1(t_2)-\widetilde B\sigma_1 D^{-1}(t_2)C(t_2) )x(t_1,t_2) + \widetilde{B}(t_2) \widetilde D^{-1}(t_2)
    \sigma_{1*}(t_2) ~y(t_1,t_2) \\
    \frac{\partial}{\partial t_2}x(t_1,t_2) = (A_2(t_2)-\widetilde B\sigma_2 D^{-1}(t_2)C(t_2) ) x(t_1,t_2) + \widetilde{B}(t_2) \widetilde D^{-1}(t_2)
    \sigma_{2*}(t_2)~y(t_1,t_2) \\
    u(t_1,t_2) = - D^{-1}(t_2) C(t_2) x(t_1,t_2) + D^{-1}(t_2) y(t_1,t_2)
    \end{array} \right.
\]
This suggests the following theorem
\begin{thm} The following collection
\[ \begin{array}{lll}
\mathfrak{V}^\times = (A_1(t_2)^\times,  ~A_2^\times(t_2), ~\widetilde B^\times(t_2), ~C^\times(t_2),~D^\times(t_2),
\widetilde D^\times(t_2); 
	\sigma_{1*}(t_2), \sigma_{2*}(t_2)\gamma_*(t_2), \sigma_1(t_2), \sigma_2(t_2), \gamma(t_2) ; \\
~~~~~~~~~~~~~~~~\mathcal H, \mathcal E_*,\mathcal{E}, \mathcal{\widetilde E_*}, \mathcal{\widetilde E})
\end{array} \]
where
\[ \begin{array}{lll}
A_1(t_2)^\times = A_1(t_2)-\widetilde B\sigma_1 D^{-1}(t_2)C(t_2), ~A_2^\times(t_2) = A_2(t_2)-\widetilde B\sigma_2 D^{-1}(t_2)C(t_2), \\
\widetilde B^\times(t_2) = \widetilde{B}(t_2) \widetilde D^{-1}(t_2), ~C^\times(t_2) = -D^{-1}(t_2)C(t_2), \\
D^\times(t_2) = D^{-1}(t_2), \widetilde D^\times(t_2) = \widetilde D^{-1}(t_2)
\end{array} \]
is a vessel with transfer function $S^\times(\lambda,t_2)$ equal to the inverse of the transfer function $S(\lambda,t_2)$ of the
vessel $\mathfrak{V}$.
\end{thm}
\textbf{Proof:} We have to show that all the vessel conditions hold. Let us omite the $t_2$ dependence of all the operators

\begin{itemize}
\item Lax equation: $\frac{d}{dt_2} A_1^\times = A_1^\times A_2^\times(t_2) - A_2^\times(t_2) A_1^\times(t_2)$.
\[ \begin{array}{llllll}
\frac{d}{dt_2} A_1^\times = \frac{d}{dt_2}[A_1-\widetilde B\sigma_1 D^{-1}C] = 
\frac{d}{dt_2}A_1 - \frac{d}{dt_2}[\widetilde B\sigma_1 D^{-1}C] 
= \frac{d}{dt_2}A_1 - \frac{d}{dt_2}[\widetilde B\sigma_1] D^{-1}C - 
\widetilde B\sigma_1\frac{d}{dt_2}[D^{-1}C]
\end{array} \]
On the other hand,
\[ \begin{array}{llllll}
A_1^\times A_2^\times(t_2) - A_2^\times(t_2) A_1^\times(t_2) = (A_1-\widetilde B\sigma_1 D^{-1}C)(A_2-\widetilde B\sigma_2 D^{-1}C) -
(A_2-\widetilde B\sigma_2 D^{-1}C)(A_1-\widetilde B\sigma_1 D^{-1}C) = \\
= A_1 A_2 - A_2 A_1 + \\
~~~~~+ [-A_1\widetilde B \sigma_2 + A_2 \widetilde B\sigma_1 - \widetilde B \sigma_2 D^{-1} C \widetilde B\sigma_1] D^{-1} C + \\
~~~~~~~~~+ \widetilde B \sigma_1 [-D^{-1} C A_2 + \sigma_1^{-1} \sigma_2 D^{-1} C A_1 + D^{-1} C \widetilde B\sigma_2 D^{-1} C]
\end{array} \]
Using now linkage conditions (\ref{eq:LinkCond}) for the original vessel $\mathfrak{DV} $ the result follows. 

\item Input vessel condition: $\frac{d}{dt_2} \big(\widetilde{B}^\times \sigma_{1*}\big) - A^\times_2 \widetilde{B}^\times \sigma_{1*} + 
A_1^\times \widetilde{B}^\times \sigma_{2*} + \widetilde{B}^\times \gamma_* = 0$. We first simplify the expression with derivative:
\[ \frac{d}{dt_2} \big(\widetilde{B}^\times \sigma_{1*}\big) =
\frac{d}{dt_2} \big( \widetilde{B} \widetilde D^{-1} \sigma_{1*}\big) =
\frac{d}{dt_2} \big( \widetilde{B} \sigma_1 D^{-1}\big) =\frac{d}{dt_2} \big( ~\widetilde{B} \sigma_1 \big) D^{-1} +
\widetilde{B} \sigma_1 \frac{d}{dt_2} \big(D^{-1}\big)
\]
and the other elements are
\[ \begin{array}{lll}
-A^\times_2 \widetilde{B}^\times \sigma_{1*} + A_1^\times \widetilde{B}^\times \sigma_{2*} + \widetilde{B}^\times \gamma_* = \\
= - [A_2-\widetilde B\sigma_2 D^{-1}C] \widetilde{B} \widetilde D^{-1} \sigma_{1*} +
 [A_1-\widetilde B\sigma_1 D^{-1}C] \widetilde{B} \widetilde D^{-1} \sigma_{2*} + \widetilde{B} \widetilde D^{-1}\gamma_* = \\
= [-A_2\widetilde B \sigma_1 + A_1 \widetilde B \sigma_2 - \widetilde B \widetilde D^{-1} \sigma_{2*} C \widetilde B \sigma_1 -
 \widetilde B \widetilde D^{-1} \sigma_{1*} C \widetilde B \sigma_2] D^{-1} + \widetilde B \widetilde D^{-1} \gamma_*
\end{array} \]
and using now the linkage condition we obtain that their sum is zero.
\item Output vessel condition: $\frac{d}{dt_2} \big(\sigma_1 C^\times \big) + \sigma_1 C^\times A^\times_2 +\sigma_2 C^\times A^\times_1+ \gamma C^\times = 0$ is similar to the input vessel condition.
\item Linkage conditions: 
\[  \begin{array}{ll}
\sigma_1 D^\times = \widetilde D^\times \sigma_{1*}, ~~~~~~ \sigma_2 D^\times = \widetilde D^\times \sigma_{2*},  \\
\widetilde D^\times \gamma_* = \sigma_2 C^\times \widetilde B^\times \sigma_{1*}  - \sigma_1 C^\times \widetilde B^\times \sigma_{2*} + 
\sigma_1 \frac{d}{dt_2} D^\times + \gamma D^\times.
\end{array}
\]
This is an immediate result of the linkage condition for the original vessel, rearranging the elements and multiplying by inverse of $D$.
\end{itemize}
\qed

The analogue of the last theorem for vessels in the integral form can be obtained by integrating
the corresponding formulas of differential vessels. The result is as follows
\begin{thm} \label{thm:InvVessel} Given the integral vessel $\mathfrak{IV}$
\[ \begin{array}{ll}
\mathfrak{IV} = (A_1(t_2), F(t_2,t_2^0),\widetilde B(t_2), C(t_2), D(t_2),\widetilde D(t_2); \\
~~~~~~~~~~~~~~~~~~~~~~~~~~~~~~~~~~~~~~~~~~~~~~~~\sigma_1(t_2), \sigma_2(t_2), \gamma(t_2), \sigma_{1*}(t_2), \sigma_{2*}(t_2)\gamma_*(t_2);
\mathcal H_{t_2}, \mathcal{E},\mathcal E_*, \mathcal{\widetilde E},\mathcal{\widetilde E_*})
\end{array} \]
the following collection is a vessel (called inverse)
\[ \begin{array}{ll}
\mathfrak{IV}^\times = (A_1^\times(t_2), F^\times(t_2,t_2^0),~\widetilde{B}(t_2) \widetilde D^{-1}(t_2), -D^{-1}(t_2)C(t_2), D^{-1}(t_2),
\widetilde D^{-1}(t_2); \\
~~~~~~~~~~~~~~~~~~~~~~~~~~~~~~~~~~~~~~~~~~~~~~~~\sigma_{1*}(t_2), \sigma_{2*}(t_2)\gamma_*(t_2), \sigma_1(t_2), \sigma_2(t_2), \gamma(t_2) ;
        \mathcal H, \mathcal E_*,\mathcal{E}, \mathcal{\widetilde E_*}, \mathcal{\widetilde E}),
\end{array} \]
where
\[ \begin{array}{llllll}
A_1(t_2)^\times = A_1(t_2)-\widetilde B\sigma_1 D^{-1}(t_2)C(t_2), \\
F^\times(t_2,t_2^0) = F(t_2,t_2^0) -
\int_{t_2^0}^{t_2} F(t_2,t_2^1) \tilde B(t_2^1) \sigma_2(t_2^1) D^{-1}(t_2^1) C(t_2^1) F(t_2^1,t_2^0) \, dt_2^1 - \ldots - \\
\int_{t_2^0}^{t_2} \int_{t_2^0}^{t_2^1} \cdots \int_{t_2^0}^{t_2^{n-1}}
F(t_2,t_2^1) \tilde B(t_2^1) \sigma_2(t_2^1) D^{-1}(t_2^1) C(t_2^1)
F(t_2^1,t_2^2) \tilde B(t_2^2) \sigma_2(t_2^2) D^{-1}(t_2^2) C(t_2^2) \times \\
~~~~~~~~~~F(t_2^2,t_2^3) \cdots F(t_2^{n-1},t_2^n) \tilde B(t_2^n) \sigma_2(t_2^n) D^{-1}(t_2^n) C(t_2^n)
F(t_2^n,t_2^0) \, dt_2^n \, dt_2^{n-1} \, \cdots \, dt_2^1 - \ldots.
\end{array} \]
\end{thm}
\textbf{Proof:} Notice that the formula for $F(t_2,t_2^0)$ is just Peano-Baker formula for $F(t_2,t_2^0)$
generalized from the differential equation (holding for the differential vessel)
\[ \frac{d}{dt_2} F^\times(t_2,t_2^0) = A_2^\times(t_2) F^\times(t_2,t_2^0) =
-[A_2(t_2)-\widetilde B\sigma_2 D^{-1}(t_2)C(t_2)] F^\times(t_2,t_2^0)
\]
\qed

\subsection{Projection, compression and cascade decomposition of systems}
Following the construction of cascade connection, it is natural to ask whether the reverse construction exists.
One of the main ingredients of this construction is that the state space $\mathcal H_{t_2}$ is decomposed into two
subspaces $\mathcal H = \mathcal H'_{t_2}\oplus\mathcal H''_{t_2}$, which are invariant for the following operators:
\[ \begin{array}{lll}
	A_1(t_2) \mathcal H'_{t_2} \subseteq \mathcal H'_{t_2}, ~~F(t_2,t_2^0) \mathcal H'_{t_2^0} = \mathcal H'_{t_2}, ~~\forall t_2,t_2^0 \\
	A_1^\times(t_2) \mathcal H''_{t_2} \subseteq \mathcal H''_{t_2}, ~~F^\times(t_2,t_2^0) \mathcal H''_{t_2^0} = \mathcal H''_{t_2}, ~~\forall t_2,t_2^0
\end{array} \]
In the differential case this means that $\mathcal H'_{t_2}$ is invariant under $A_1(t_2), A_2(t_2)$ for all $t_2$ and that
$\mathcal H''_{t_2}$ is invariant under $A_1^\times(t_2), A_2^\times(t_2)$ for all $t_2$.
\begin{defn} Given a vessel $\mathfrak{IV}$ subspaces $\mathcal G_{t_2} \subseteq \mathcal H_{t_2}$ are called \textbf{invariant}
if for all $t_2$ holds
\[ A_1(t_2) \mathcal G_{t_2} \subseteq \mathcal G_{t_2}, ~~F(t_2,t_2^0) \mathcal G_{t_2^0} = \mathcal G_{t_2}.
\]
Subspaces $\mathcal G^\times_{t_2} \subseteq \mathcal H_{t_2}$ are called \textbf{co-invariant} if their complements
$\mathcal H_{t_2} \dot\ominus \mathcal G^\times_{t_2}$ are invariant and if these subspaces are invariant for the operators $A_1^\times(t_2)$ and
$F^\times(t_2,t_2^0)$ (defined in theorem \ref{thm:InvVessel}), i.e., if for all $t_2$ holds
\[ A_1^\times(t_2) \mathcal G^\times_{t_2} \subseteq \mathcal G^\times_{t_2}, ~~F^\times(t_2,t_2^0) \mathcal G^\times_{t_2^0} = \mathcal G^\times_{t_2}.
\]
\end{defn}
The classical condition for a cascade decomposition and a factorization of the transfer function \cite{bib:BGK,bib:BGKD,bib:S} uses these
two notions of invariance. We present an analogue of the corresponding theorems.

Assume that we are given an overdetermined $2D$ system, $t_1$ invariant (\ref{eq:system}) with the vessel
\[ \mathfrak{IV} = (A_1(t_2), F(t_2,t_2^0),\widetilde B(t_2),C(t_2), D(t_2),\widetilde D(t_2);
    \sigma_1(t_2), \sigma_2(t_2), \gamma(t_2),  \sigma_{1*}(t_2), \sigma_{2*}(t_2)\gamma_*(t_2);
        \mathcal H_{t_2}, \mathcal{E},\mathcal E_*, \mathcal{\widetilde E},\mathcal{\widetilde E_*}).
\]
Suppose also that we are given subspaces $\mathcal G_{t_2} \subseteq \mathcal H_{t_2}$ that are invariant.
Then it is possible to define a  \textit{projection} of the vessel $\mathfrak{IV}$ onto the invariant subspaces $\mathcal G_{t_2}$
as follows
\begin{defn} \label{def:Proj}
Projection of the vessel $\mathfrak{IV}$ on the invariant subspaces $\mathcal G_{t_2}$ is a collection
\[ \mathfrak{IV'} = (A_1'(t_2), F'(t_2,t_2^0),\widetilde B'(t_2),C'(t_2), D(t_2),\widetilde D(t_2); \sigma_1(t_2), \sigma_2(t_2), \gamma(t_2),
\sigma'_{1*}(t_2), \sigma'_{2*}(t_2)\gamma'_*(t_2); \mathcal G_{t_2}, \mathcal{E},\mathcal E_*, \mathcal{\widetilde E}',\mathcal{\widetilde E_*}').
\]
where denoting by $P_{\mathcal G_{t_2}}$ - projection on $\mathcal G_{t_2}$, the operators in $\mathfrak{IV'}$ are
\[ \begin{array}{lllll}
F'(t_2,t_2^0) = F(t_2,t_2^0) P_{G(t_2^0)}, \\
A_1'(t_2) = F'(t_2,t_2^0) A_1(t_2^0) F'(t_2^0,t_2), \\
B'(t_2) = P_{\mathcal G_{t_2}} B(t_2) , \\
C'(t_2) = C(t_2) P_{\mathcal G_{t_2}},
\end{array} \]
and $\sigma'_{1*}(t_2), \sigma'_{2*}(t_2), \gamma'_*(t_2)$ are taken so that the linkage conditions (\ref{eq:LinkCond}) are satisfied
\[ \sigma'_{1*} D= \widetilde D \sigma_1, ~~ \sigma'_{2*} D = \widetilde D \sigma_2, 
~~\widetilde D \gamma = \sigma'_{2*} C' \widetilde B' \sigma_1  - \sigma'_{1*} C' \widetilde B' \sigma_2 + \sigma'_{1*} \frac{d}{dt_2}D + \gamma'_* D.
\]
\end{defn}
Let $\mathcal G^\times_{t_2}$ be co-invariant subspaces. Then we define \textit{compression} of the vessel $\mathfrak{IV}$ 
onto the the co-invariant subspaces $\mathcal G^\times_{t_2}$ as follows
\begin{defn} \label{def:Compr}
Compression of the vessel $\mathfrak{IV}$ on the co-invariant subspaces $\mathcal G^\times_{t_2}$ is a collection
\[ \begin{array}{lll}
	\mathfrak{IV^\times} = (A_1^\times(t_2), F^\times(t_2,t_2^0),	\widetilde B^\times(t_2),C^\times(t_2), D(t_2),\widetilde D(t_2); \\
	\sigma''_1(t_2), \sigma''_2(t_2), \gamma''(t_2), \sigma_{1*}(t_2), \sigma_{2*}(t_2)\gamma_*(t_2);
	\mathcal G^\times_{t_2}, \mathcal{E}'',\mathcal E''_*, \mathcal{\widetilde E},\mathcal{\widetilde E_*}).
\end{array} \]
where $A_1^\times(t_2), F^\times(t_2,t_2^0)$ are defined as in theorem \ref{thm:InvVessel} and denoting
by $P^\times_{\mathcal G^\times_{t_2}}$ - projection on $\mathcal G^\times_{t_2}$
\[ \begin{array}{lllll}
B^\times(t_2) = P^\times_{\mathcal G_{t_2}} B(t_2) , ~~C''(t_2) = C(t_2) P^\times_{\mathcal G_{t_2}},
\end{array} \]
the operators $\sigma''_1(t_2), \sigma''_2(t_2), \gamma''(t_2)$ are taken so that the linkage conditions (\ref{eq:LinkCond}) are satisfied
\[\sigma_{1*} D = \widetilde D \sigma''_1 , ~~ \sigma_{2*} D = \widetilde D \sigma''_2, 
~~\widetilde D \gamma'' = \sigma_{2*} C^\times \widetilde B^\times \sigma''_1  - \sigma_{1*} C^\times \widetilde B^\times \sigma''_2 + \sigma_{1*} \frac{d}{dt_2}D + \gamma_* D
\]
\end{defn}
\begin{lemma} $\mathfrak{IV'}, \mathfrak{IV}^\times$ are vessels.
\end{lemma}
\textbf{Proof:} We will show that $\mathfrak{IV'}$ is a vessel. For $\mathfrak{IV''}$ the proof is essentially the same.

First we show that $F'(s,t)$ is an evolution group. $F'(s,s) = F(s,s) P_{\mathcal G(s)} = I_{\mathcal G(t_2^0)}$,
\[ F'(s,t) F'(t,y) = F(s,t)P_{\mathcal G(t)} F(t,y) P_{\mathcal G(y)} = F(s,t) F(t,y) P_{\mathcal G(y)} = F(s,y) P_{\mathcal G(y)} = F'(t,y).
\]
Then we have to show that all vessel conditions are satisfied.
\[ A_1'(t_2) F'(t_2,t_2^0) = F'(t_2,t_2^0) A_1(t_2^0) F'(t_2,t_2^0) F'(t_2,t_2^0) = F(t_2,t_2^0) A_1(t_2^0) P_{\mathcal G(t_2^0)} =
F'(t_2,t_2^0) A'_1(t_2^0),
\]
which means that the Lax equation holds. In order to check the input vessel condition notice that
\[ F' (t_2^0, t_2) \widetilde B'(t_2) = P_{\mathcal G(t_2^0)} F (t_2^0, t_2) P_{\mathcal G_{t_2}} \widetilde B(t_2) =
P_{G(t_2^0)} F (t_2^0, t_2) \widetilde B(t_2),
\]
and in the same manner
\[ \begin{array}{lll}
F'(t_2^0, t_2) A'_1(t_2) \widetilde B'(t_2) = F'(t_2^0,t_2) F'(t_2,t_2^0) A_1(t_2^0) F'(t_2^0,t_2) P_{\mathcal G_{t_2}} B(t_2) = \\
= A_1(t_2^0) P_{\mathcal G(t_2^0)} F(t_2^0, t_2) \widetilde B(t_2) = P_{\mathcal G(t_2^0)} F(t_2^0, t_2) A_1(t_2) \widetilde B(t_2).
\end{array} \]
So, the input vessel condition is
\[ \begin{array}{lll}
\frac{d}{dt_2} (F' (t_2^0, t_2) \widetilde B'(t_2) \sigma_1(t_2)) + F'(t_2^0, t_2) A'_1(t_2) \widetilde B'(t_2) \sigma_2(t_2)
   + F'(t_2^0, t_2) \widetilde B'(t_2) \gamma(t_2) =  \\
=  P_{G(t_2^0)}[\frac{d}{dt_2}( F (t_2^0, t_2) \widetilde B(t_2) \sigma_1(t_2)) + F(t_2^0, t_2) A_1(t_2) \widetilde B(t_2) \sigma_2(t_2)
   + F(t_2^0, t_2) \widetilde B(t_2) \gamma(t_2)] = P_{G(t_2^0)} 0 = 0.
\end{array} \]
The output vessel condition is a consequence of other conditions by considering the differential equation
\[ \frac{\partial}{\partial t_2} S'(\lambda,t_2) =
\sigma_{1*}^{-1}(t_2) (\lambda \sigma_{2*}(t_2) + \gamma'_{*}(t_2)) S'(\lambda,t_2) -
S'(\lambda,t_2) \sigma_1^{-1}(t_2) (\lambda \sigma_2(t_2) + \gamma(t_2))
\]
for the projected transfer function
\begin{equation} \label{eq:DifEqS'}
\begin{array}{llll}
S'(\lambda,t_2) & = D(t_2) + C'(t_2)(\lambda I - A_1'(t_2))^{-1} \widetilde B'(t_2) \sigma_1 = \\
	& = D(t_2) + C(t_2) F(t_2,t_2^0)P_{t_2^0}(\lambda I - A_1(t_2^0))^{-1} P_{t_2^0} F(t_2^0,t_2)\widetilde B(t_2) \sigma_1.
\end{array} \end{equation}
Let us first evaluate the derivative
\[ \begin{array}{llll}
\frac{\partial}{\partial t_2} S'(\lambda,t_2) =
\frac{\partial}{\partial t_2} D + \frac{\partial}{\partial t_2} [C(t_2) F(t_2,t_2^0)P_{t_2^0}](\lambda I - A_1(t_2^0))^{-1} P_{t_2^0} F(t_2^0,t_2)\widetilde B(t_2) \sigma_1 + \\
~~~~~ + C(t_2) F(t_2,t_2^0)P_{t_2^0}(\lambda I - A_1(t_2^0))^{-1} P_{t_2^0}\frac{\partial}{\partial t_2}[F(t_2^0,t_2)\widetilde B(t_2) \sigma_1] = \\
= \frac{\partial}{\partial t_2} D + \sigma_{1*}^{-1}[\sigma_{2*}C(t_2) F(t_2,t_2^0) A_1(t_2^0)P_{t_2^0} +
\gamma'C(t_2) F(t_2,t_2^0)P_{t_2^0}](\lambda I - A_1(t_2^0))^{-1} P_{t_2^0} F(t_2^0,t_2)\widetilde B(t_2) \sigma_1 - \\
~~~~~ -C(t_2) F(t_2,t_2^0)P_{t_2^0}(\lambda I - A_1(t_2^0))^{-1} P_{t_2^0}[A_1(t_2^0)F(t_2^0,t_2)\widetilde B(t_2) \sigma_2 + F(t_2^0,t_2)\widetilde B(t_2) \gamma] = \\
= \frac{\partial}{\partial t_2} D + \sigma_{1*}^{-1}[\sigma_{2*}\lambda + \gamma']C(t_2) F(t_2,t_2^0)P_{t_2^0}](\lambda I - A_1(t_2^0))^{-1} P_{t_2^0} F(t_2^0,t_2)\widetilde B(t_2) \sigma_1 - \\
~~~~~ -C(t_2) F(t_2,t_2^0)P_{t_2^0}(\lambda I - A_1(t_2^0))^{-1} P_{t_2^0}A_1(t_2^0)F(t_2^0,t_2)\widetilde B(t_2) [\lambda \sigma_2 + \gamma] - \\
~~~~~~~~~~- \sigma_{1*}^{-1}\sigma_{2*}C(t_2) F(t_2,t_2^0)P_{t_2^0}F(t_2^0,t_2)\widetilde B(t_2) \sigma_1 +
 C(t_2) F(t_2,t_2^0)P_{t_2^0}F(t_2^0,t_2)\widetilde B(t_2)\sigma_2 = \\
= \frac{\partial}{\partial t_2} D + \sigma_{1*}^{-1}[\sigma_{2*}\lambda + \gamma'](-D + S'(\lambda,t_2)) 
 -(-D + S'(\lambda,t_2))\sigma_1^{-1} [\lambda \sigma_2 + \gamma] - \\
~~~~~~~~~~- \sigma_{1*}^{-1}\sigma_{2*}C(t_2) F(t_2,t_2^0)P_{t_2^0}F(t_2^0,t_2)\widetilde B(t_2) \sigma_1 +
 C(t_2) F(t_2,t_2^0)P_{t_2^0}F(t_2^0,t_2)\widetilde B(t_2)\sigma_2 = \\
= \frac{\partial}{\partial t_2} D + \sigma_{1*}^{-1}[\sigma_{2*}\lambda + \gamma'](-D + S'(\lambda,t_2)) 
 -(-D + S'(\lambda,t_2))\sigma_1^{-1} [\lambda \sigma_2 + \gamma] - \\
~~~~~~~~~~- \sigma_{1*}^{-1}\sigma_{2*}C'(t_2) P_{t_2^0}\widetilde B'(t_2) \sigma_1 +
 C'(t_2) P_{t_2^0}\widetilde B'(t_2)\sigma_2.
\end{array} \]
If we plug this expression into (\ref{eq:DifEqS'}), we shall obtain
\[ \begin{array}{llll}
\frac{\partial}{\partial t_2} D - \sigma_{1*}^{-1}[\sigma_{2*}\lambda + \gamma']D + D \sigma_1^{-1} [\lambda \sigma_2 + \gamma]- \\
~~~~~~~~~~- \sigma_{1*}^{-1}\sigma_{2*}C'(t_2) P_{t_2^0}\widetilde B'(t_2) \sigma_1 +
 C'(t_2) P_{t_2^0}\widetilde B'(t_2)\sigma_2 = 0
\end{array} \]
\qed

It is also possible to perform a compression of the vessel on a semi-invariant subspace. Let us first define it. Suppose that 
$\mathcal G_{t_2}$ is an invariant subspace of the vessel $\mathfrak{IV}$. Suppose also that there is a co-invariant within $\mathcal G_{t_2}$
subspace $\mathcal G'_{t_2} \subseteq \mathcal G_{t_2}$, that is the subspace $\mathcal G_{t_2} \dot- \mathcal G'_{t_2}$ is co-invariant.
In this case we call $\mathcal G'_{t_2}$ a \textit{semi-invariant} subspace of $\mathfrak{IV}$. Then performing projection of the vessel
$\mathfrak{IV}$ on the invariant subspace $\mathcal G_{t_2}$, we shall obtain a new vessel $\mathfrak{IV}'$. Performing further compression on
the co-invariant subspace $\mathcal G'_{t_2}$ we shall obtain the desired vessel $\mathfrak{IV}''$.

\noindent\textbf{Remarks: 1.} Notice that invertibility of $D(t_2)$ is not essential but suffices in order to determine
all the relevant data. For example, in the case of projection on an invariant space $\mathcal G_{t_2}$ one obtains
\[ \sigma'_{1*} = \widetilde D \sigma_1 D^{-1}, ~~ \sigma'_{2*} = \widetilde D \sigma_2 D^{-1},
~~\gamma'_* = \widetilde D \gamma D^{-1} - [ \sigma'_{2*} C' \widetilde B' \sigma_1  - \sigma'_{1*} C' \widetilde B' \sigma_2 + \sigma'_{1*} \frac{d}{dt_2}D] D^{-1} 
\]
and consequently they are uniquely determined.

\noindent\textbf{2.} It is enough to have a subspace $\mathcal G_{t_2^0} \subseteq \mathcal H_{t_2^0}$ that is invariant
under $A_1(t_2^0)$ with a complementary subspace $\mathcal G^\times_{t_2^0} \subseteq \mathcal H_{t_2^0}$ that is invariant under 
$A_1^\times(t_2^0)$ and such that
\[ F(t_2,t_2^0) \mathcal G_{t_2^0}, ~~F(t_2,t_2^0) \mathcal G^\times_{t_2^0}
\]
are complementary for all $t_2$.

\noindent\textbf{3.} For the conservative case (to be studied later) $A^\times_1(t_2) = A_1^*(t_2)$, so the existence of a complementary 
invariant subspace is automatic, just like in the classical case.

As a result of all these consideration we deduce a theorem of cascade decomposition of vessels. Suppose that we are given 
an invariant subspace $\mathcal G_{t_2}$, which is at the same time co-invariant. Then it is possible to produce a projection
on $\mathcal G_{t_2}$ and compression on its complement. The point is that the obtained in such a way vessel can be cascade connected to give
the initial one. This is precisely the content of the next theorem
\begin{thm} Suppose that we are given a vessel $\mathfrak{IV}$ of the form (\ref{eq:CascVessel}) with invariant subspace
$\mathcal G_{t_2}$, and co-invariant subspace $\mathcal G^\times_{t_2}$, then the projection on $\mathcal G_{t_2}$ produces a vessel
$\mathfrak{IV}'$ and compression to $\mathcal G^\times_{t_2}$ produces $\mathfrak{IV}^\times$.
Moreover, it is possible to cascadly connect these two vessels and to obtain the original one $\mathfrak{IV}$.
\end{thm}
\section{\label{sec:Kalman}Kalman decomposition of Vessels}
The notions of approximate controllability and observability
allows building of minimal systems for which there is a very
good classification theory. It turns out that there are possible
other parts of a system (vessel) that are
non-approximately controllable or are not observable. Let us
denote the system (\ref{eq:system}) (for fixed $\sigma_1(t_2),
\sigma_2(t_2), \gamma(t_2), \sigma_{1*}(t_2), \sigma_{2*}(t_2),
\gamma_*(t_2)$)
\[
     I\Sigma: \left\{ \begin{array}{lll}
    \frac{\partial}{\partial t_1}x(t_1,t_2) = A_1(t_2) ~x(t_1,t_2) + \widetilde B(t_2) \sigma_1(t_2) ~u(t_1,t_2) \\[5pt]
    x(t_1, t_2) = F (t_2,t_2^0) x(t_1, t_2^0) + \int\limits_{t_2^0}^{t_2} F(t_2, s) \widetilde B(s) \sigma_2(s) u(t_1, s)ds \\[5pt]
    y(t_1,t_2) = C(t_2)~ x(t_1,t_2) + D(t_2) u(t_1,t_2).
  \end{array} \right.
\]
by
\[ I\Sigma = [A_1(t_2), F(t_2,t_2^0), \widetilde B(t_2), C(t_2), D(t_2)].
\]
\begin{thm}
Let $I\Sigma$ be a system defined in (\ref{eq:system}) with inner state space $\mathcal H_{t_2}$ for each $t_2$.
Then there exists an orthogonal (for each $t_2$) decomposition of the state space
\[ \mathcal H_{t_2} = \mathcal H_{t_2}^{c\bar o} \oplus \mathcal H_{t_2}^{co} \oplus \mathcal H_{t_2}^{\Bar{co}} \oplus \mathcal H_{t_2}^{\bar co}
\]
so that with respect to this decomposition the system has the following decomposition of its operators
\[ \begin{array}{lllll}
A_1(t_2) = \bbmatrix{A_1^{c\bar o}(t_2) & A_{12}(t_2) & A_{13}(t_2) & A_{14}(t_2) \\ 0 & A_1^{co}(t_2) & A_{23}(t_2) & A_{24}(t_2) \\
                             0 & 0 & A_1^{\Bar{co}}(t_2) & A_{34}(t_2) \\   0 & 0 & 0 & A_1^{\bar co}(t_2)}; \\
    \widetilde B(t_2) = \bbmatrix{\widetilde B^{c\bar o} \\ \widetilde B^{co} \\ 0 \\ 0 } \\
    C(t_2) = \bbmatrix{0 & C^{co}(t_2) & 0 & C^{\bar c o}(t_2)} \\
    F(t_2,t_2^0) = \bbmatrix{F^{c\bar o}(t_2,t_2^0) \\ F^{co}(t_2,t_2^0) \\ F^{\Bar{co}}(t_2,t_2^0) \\  F^{\bar co}(t_2,t_2^0) }
\end{array} \]
where the subsystem, defined by
\[ I\Sigma^c = [ \bbmatrix{A_1^{c\bar o}(t_2) & A_{12}(t_2) \\ 0 & A_1^{co}(t_2)}, \bbmatrix{F^{c\bar o}(t_2,t_2^0) \\ F^{co}(t_2,t_2^0)},
         \bbmatrix{\widetilde B^{c\bar o} \\ \widetilde B^{co}}, \bbmatrix{0 & C^{co}(t_2)}, D ]
\]
is approximately controllable, the system
\[ I\Sigma^o = [\bbmatrix{A_1^{co}(t_2) & A_{24}(t_2) \\ 0 & A_1^{\bar co}(t_2)}, \bbmatrix{F^{co}(t_2,t_2^0) \\ F^{\bar co}(t_2,t_2^0)},
         \bbmatrix{\widetilde B^{co} \\ 0}, \bbmatrix{C^{co}(t_2) & C^{\bar co}(t_2)}, D]
\]
is observable, and the system
\[ I\Sigma^{co} = [A_1^{co}, F^{co}(t_2,t_2^0), \widetilde B^{co}(t_2), C^{co}(t_2), D]
\]
is minimal (i.e., approximately controllable and observable). Moreover, the transfer functions of all the systems are equal:
\[ S(\lambda,t_2) = S^c(\lambda,t_2) = S^o(\lambda,t_2) = S^{co}(\lambda,t_2).
\]
\end{thm}
\textbf{Proof:} Denote for each $t_2$
\[ \begin{array}{lll}
\mathcal G_c(t_2) = \operatorname{\overline{span}} \{ \operatorname{Im} A_1^j(t_2) B(t_2) \mid j=0,1,2,\ldots \}, &
 \mathcal G_{\bar c}(t_2) = \mathcal H_{t_2} \ominus \mathcal G_c(t_2) \\
\mathcal G_o(t_2) = \operatorname{\overline{span}} \{ \operatorname{Im} A_1^{*j}(t_2) C^*(t_2) \mid j=0,1,2,\ldots \}, &
\mathcal G_{\bar o}(t_2) = \mathcal H_{t_2} \ominus \mathcal G_o(t_2) \\
\end{array} \]
Now write each operator in the system $I\Sigma$
as an operator with respect to $\mathcal G_c(t_2)$ and $\mathcal G_{\bar c}(t_2)$. Since $\mathcal G_c(t_2)$ is
$A_1(t_2)$ invariant, the system $I\Sigma = [A_1(t_2), F(t_2,t_2^0), \widetilde B(t_2), C(t_2), D(t_2)]$ is of the following form
\[ [\bbmatrix{A_1^c(t_2) & A_{12}(t_2) \\ 0 & A_1^{\bar c}(t_2)}, \bbmatrix{ F^c(t_2,t_2^0)\\F^{\bar c}(t_2,t_2^0)},
\bbmatrix{\widetilde B^c(t_2)\\ \widetilde B^{\bar c}(t_2)}, \bbmatrix{C^c(t_2) & C^{\bar c}(t_2)}, D(t_2)]
\]
and clearly
\[ C^c(t_2) A^c(t_2) B^c(t_2) = C(t_2) A(t_2) B(t_2).
\]
Thus the original system $I\Sigma$ has the same transfer function as
\[ I\Sigma^c = [A_1^c(t_2), F^c(t_2,t_2^0), \widetilde B^c(t_2), C^c(t_2), D(t_2)]
\]
and this system is approximately controllable. The same process works on the given system (\ref{eq:system}) with observability
($\mathcal G_{\bar o}$ rather than $\mathcal G_o$ is invariant under $A_1(t_2)$) to give a ’’smaller'' observable system
\[ I\Sigma^o = [A_1^o(t_2), F^o(t_2,t_2^0), \widetilde B^o(t_2), C^o(t_2), D(t_2)].
\]
If one combines these two processes, one gets (by first decomposing the system $I\Sigma$ into controllable and uncontrollable parts and
then decomposing these systems into observable and unobservable parts) the desired decomposition. \qed

\section{\label{sec:AnalFunc}Analytic functions as transfer functions}
The aim of this section is to show that any function in our class
$\boldsymbol{\mathcal I}$ can be realized, i.e., presented as a
transfer function of a certain vessel.

\subsection{Realization theorem for arbitrary analytic functions in $\boldsymbol{\mathcal I}$}
So, suppose that we are given a function $S(\lambda,t_2)\in \boldsymbol{\mathcal I}$. Our first aim is to realize this function, i.e. to show that this class is realizable. In order to do it, we realize $S(\lambda,t_2)$ for a fixed $t_2^0$ \cite{bib:Helton} as
\[ S(\lambda,t_2^0) = D_0 + C_0 (\lambda I - A_1)^{-1} B_0 \sigma_1(t_2^0).
\]
Then the following theorem holds
\begin{thm} \label{tm:NonConsRealiz}
Suppose that $S(\lambda, t_2) \in \boldsymbol{\mathcal I}$. Then there exists vessel $\mathfrak{DV}$ in the differential form
\[ \mathfrak{DV} = (A_1(t_2), A_2(t_2), B(t_2), C(t_2), D(t_2), \widetilde D;
        \sigma_1, \sigma_2, \gamma, \sigma_{1*}, \sigma_{2*}, \gamma_*;
                \mathcal{H}, \mathcal{E}, \mathcal{E}_*, \widetilde{\mathcal E}, \widetilde{\mathcal E}_*),
\]
with this transfer function and for which
\begin{eqnarray}
C(t_2) = \oint\limits_{Spec A_1} \Phi_*(\lambda,t_2,t_2^0) C_0(\lambda I - A_1)^{-1} d\lambda \\
B(t_2) = \oint\limits_{Spec A_1} (\lambda I - A_1)^{-1} B_0 \sigma_1(t_2) \Phi^* (-\bar\lambda,t_2,t_2^0) d\lambda \\
D(t_2) = S(\infty, t_2)
\end{eqnarray}
and
\[ S(\lambda, t_2) = D(t_2) + C(t_2) (\lambda I - A_1)^{-1} B(t_2) \sigma_1. \]
\end{thm}
Before we consider the proof of this theorem, it is important to note that for the functions $C(t_2), B(t_2)$ the following lemma holds
\begin{lemma} \label{lm:DEForCB} $C(t_2), B(t_2)$ satisfy the following differential equations with the spectral operator parameter $A_1$:
\[ \begin{array}{llll}
\frac{d}{dt_2} C(t_2) = \sigma_{1*}^{-1}(t_2) [\sigma_{2*}(t_2) C(t_2) A_1 + \gamma_*(t_2) C(t_2)] \\
\frac{d}{dt_2} [B(t_2) \sigma_1(t_2)] = A_1 B(t_2) \sigma_2(t_2) + B(t_2) \gamma(t_2)
\end{array} \]
\end{lemma}
\textbf{Proof:} For $C(t_2)$ one obtains that
\[ \begin{array}{llll}
 \sigma_1 C(t_2)' & = \sigma_1 \frac{\partial}{\partial t_2}\oint\limits_{Spec A_1} \Phi_*(\lambda,t_2,t_2^0) C_0(\lambda I - A_1)^{-1} d\lambda = \\
        & = \sigma_1 \oint\limits_{Spec A_1} \sigma_1^{-1} (\sigma_2 \lambda + \gamma_*) \Phi_*(\lambda,t_2,t_2^0) C_0(\lambda I - A_1)^{-1} d\lambda = \\
        & = \oint\limits_{Spec A_1} \sigma_2 \Phi_*(\lambda,t_2,t_2^0) C_0 (\lambda I - A_1 + A_1) (\lambda I - A_1)^{-1} d\lambda + \\
    & ~~~~~~~~~ + \oint\limits_{Spec A_1} \gamma_* \Phi_*(\lambda,t_2,t_2^0) C_0(\lambda I - A_1)^{-1} d\lambda = \\
        & =\sigma_2  \oint\limits_{Spec A_1}  \Phi_*(\lambda,t_2,t_2^0) C_0 (\lambda I - A_1)^{-1} d\lambda A_1 +
       \gamma_* \oint\limits_{Spec A_1} \Phi_*(\lambda,t_2,t_2^0) C_0(\lambda I - A_1)^{-1} d\lambda = \\
     & = \sigma_2 C(t_2) A_1 + \gamma_* C(t_2)
\end{array}
\]
and the same proof works for $B(t_2)$. \qed

\noindent\textbf{Proof of theorem \ref{tm:NonConsRealiz}:}
Let us define a vessel:
\[ \mathfrak{DV} =
(A_1, A_2=0, B(t_2), C(t_2), D(t_2), \widetilde D = \sigma_{1*} D(t_2) \sigma_1^{-1}; \sigma_1, \sigma_2, \gamma, \sigma_{1*}, \sigma_{2*}, \gamma_*;
 \mathcal{H}, \mathcal{E}, \mathcal{E}_*, \widetilde{\mathcal E}, \widetilde{\mathcal E}_*),
\]
and show that all vessel conditions hold. Lax equation holds
\[ \frac{d}{dt_2} A_1 = 0 = A_1 A_2 - A_2 A_1 = A_1 0 - 0 A_1.
\]
The first (\ref{eq:OverDetCondIn}) and the second (\ref{eq:OverDetCondOut}) vessel conditions are exactly the contents of
lemma \ref{lm:DEForCB}.

Consider next the expression $C(t_2)(\lambda I - A_1)^{-1} B(t_2) \sigma_1(t_2)$. Using lemma (\ref{lm:DEForCB}) we obtain that
\[ \begin{array}{llllll}
\frac{\partial}{\partial t_2} \big( C(t_2)(\lambda I - A_1)^{-1} B(t_2) \sigma_1(t_2) \big) = \\
\sigma_{1*}^{-1}(\sigma_{2*} C(t_2) A_1 + \gamma_*(t_2) C(t_2))(\lambda I - A_1)^{-1} B(t_2) \sigma_1(t_2) + \\
~~~~~~~~~~~~~ + C(t_2)(\lambda I - A_1)^{-1} [A_1 B(t_2) \sigma_2(t_2) + B(t_2) \gamma(t_2)] = \\
= \sigma_{1*}^{-1}(\sigma_{2*} \lambda + \gamma_*(t_2)) C(t_2)(\lambda I - A_1)^{-1} B(t_2) \sigma_1(t_2) -
C(t_2)(\lambda I - A_1)^{-1} B(t_2) (\sigma_2(t_2) \lambda+\gamma(t_2)) - \\
~~~~~~~~~~~~~ - \sigma_{1*}^{-1}\sigma_{2*} C(t_2) B(t_2) \sigma_1(t_2) + C(t_2) B(t_2) \sigma_1(t_2) \sigma_1^{-1}(t_2) \sigma_2(t_2)
\end{array} \]
We can use here the fundamental matrices $\Phi(\lambda,t_2,t_2^0), \Phi_*(\lambda,t_2,t_2^0)$ and with the help of variation of
coefficients, we obtain that
\[ C(t_2)(\lambda I - A_1)^{-1} B(t_2) \sigma_1(t_2) = \Phi_*(\lambda,t_2,t_2^0) K(\lambda,t_2) \Phi^{-1}(\lambda,t_2,t_2^0),
\]
where
\[ \begin{array}{llll}
\frac{\partial}{\partial t_2} K(\lambda,t_2) = \\
\Phi_*^{-1}(\lambda,t_2,t_2^0) \big( - \sigma_{1*}^{-1}\sigma_{2*} C(t_2) B(t_2) \sigma_1(t_2) + C(t_2) B(t_2) \sigma_1(t_2) \sigma_1^{-1}(t_2) \sigma_2(t_2) \big) \Phi(\lambda,t_2,t_2^0)
\end{array} \]
with the initial value $K(\lambda,t_2^0) = C(t_2^0)(\lambda I - A_1)^{-1} B(t_2^0) \sigma_1(t_2^0)$. Since the fundamental matrices
are entire in $\lambda$ functions, we also obtain that
\[ K(\lambda,t_2) - K(\lambda,t_2^0) = \int_{t_2^0}^{t_2} \frac{\partial}{\partial y} K(\lambda,y) dy
\]
is an entire function of $\lambda$. Thus
\[ \begin{array}{llll}
S(\lambda,t_2) - C(t_2)(\lambda I - A_1)^{-1} B(t_2) \sigma_1(t_2) =
\Phi_*(\lambda,t_2,t_2^0) \big( S(\lambda,t_2^0) - K(\lambda,t_2) \big) \Phi^{-1}(\lambda,t_2,t_2^0) = \\
= \Phi_*(\lambda,t_2,t_2^0) \big( D(t_2^0) + C(t_2^0)(\lambda I - A_1)^{-1} B(t_2^0) \sigma_1(t_2^0) - \\
~~~~~~~~~~~~~~~~~~~~~~~~~
-C(t_2^0)(\lambda I - A_1)^{-1} B(t_2^0) \sigma_1(t_2^0) - \int_{t_2^0}^{t_2}
        \frac{\partial}{\partial y} K(\lambda,y) dy \big) \Phi^{-1}(\lambda,t_2,t_2^0) = \\
= \Phi_*(\lambda,t_2,t_2^0) \big( D(t_2^0) + \frac{\partial}{\partial y} K(\lambda,y) dy \big) \Phi^{-1}(\lambda,t_2,t_2^0),
\end{array} \]
which is an entire in $\lambda$ function too. On the other hand, when $\lambda$ tends to $\infty$, $S(\infty,t_2) = D(t_2)$ and
$\lim\limits_{\lambda\rightarrow\infty}C(t_2)(\lambda I - A_1)^{-1} B(t_2) \sigma_1(t_2) = 0$, which means that their
difference is bounded and consequently, by Liouville’s theorem for operator valued functions, applied to constant vectors
is constant. Finally,
\begin{equation} \label{eq:SequalDCAB}
 S(\lambda,t_2) = D(t_2) + C(t_2)(\lambda I - A_1)^{-1} B(t_2) \sigma_1(t_2).
\end{equation}
Once we have established all these formulas, it remains to show that the linkage conditions
(\ref{eq:LinkCond})
\[ \begin{array}{llll}
\sigma_{1*} D = \widetilde D \sigma_1, ~~ \sigma_{2*} D = \widetilde D \sigma_2 \\
\widetilde D \gamma = \sigma_{2*} C \widetilde B \sigma_1 -\sigma_{1*} C \widetilde B \sigma_2 -
\sigma_{1*} \frac{d}{dt_2} D + \gamma_* D
\end{array} \]
are satisfied. In order to do this, we use the differential equation for $S(\lambda,t_2)$ (\ref{eq:DforS})
\[
\frac{\partial}{\partial t_2} S(\lambda,t_2) = \sigma_{1*}^{-1}(t_2) (\sigma_{2*}(t_2) \lambda + \gamma_*(t_2)) S(\lambda,t_2)-
S(\lambda,t_2)\sigma_1^{-1}(t_2) (\sigma_2(t_2) \lambda + \gamma(t_2)).
\]
and substitute here the realization formula (\ref{eq:SequalDCAB}). Then
\[ \begin{array}{lll}
\frac{d}{dt_2} D(t_2) + \frac{d}{dt_2} C(t_2)(\lambda I - A_1)^{-1} B(t_2) \sigma_1(t_2) +
C(t_2)(\lambda I - A_1)^{-1} \frac{d}{dt_2} \big( B(t_2) \sigma_1(t_2) \big) = \\
= \sigma_{1*}^{-1}(t_2) (\sigma_{2*}(t_2) \lambda + \gamma_*(t_2)) \big( D(t_2) + C(t_2)(\lambda I - A_1)^{-1} B(t_2) \sigma_1(t_2)\big) - \\
~~~~ - \big( D(t_2) + C(t_2)(\lambda I - A_1)^{-1} B(t_2) \sigma_1(t_2) \big) \sigma_1^{-1}(t_2) (\sigma_2(t_2) \lambda + \gamma(t_2)).
\end{array} \]
Considering the linear in $\lambda$ part, for big values of $\lambda$ we immediately obtain that
\[ \sigma_{1*}^{-1} \sigma_{2*} D = D \sigma_1^{-1}\sigma_2,
\]
and plugging this back into the differential equation, and tending $\lambda$ to infinity,
we obtain (defining $\widetilde D = \sigma_{1*} D\sigma_1^{-1}$) that
\[ \widetilde D \gamma = \sigma_{2*} C \widetilde B \sigma_1 -\sigma_{1*} C \widetilde B \sigma_2 -
\sigma_{1*} \frac{d}{dt_2} D + \gamma_* D,
\]
which finishes the proof. \qed

Notice that there are no assumptions on the dimensions of $\mathcal E,\mathcal E_*$. In the next section we consider the finite dimensional
case $\dim\mathcal H_{t_2}<\infty$.

\subsection{Realization theorem for matrix functions in $\boldsymbol{\mathcal I}$}
In this section we want to further investigate the formulas for $C(t_2), B(t_2)$ arising in the realization of
$S(\lambda,t_2)\in\mathrm K$ for the matrix case. So, $S(\lambda, t_2)$ maps solutions of the input
ODE (\ref{eq:InCC}) with the spectral parameter $\lambda$ to solutions of the output ODE (\ref{eq:OutCC}) with the same spectral parameter.
By Jordan theorem any constant matrix is similar to its unique Jordan form, and consequently, we are going to concentrate on this special case.
More explicit realization, based on the theorem \ref{tm:NonConsRealiz} is achieved in the next theorem.

So, suppose that there is only one eigenvalue for $A_1$ of multiplicity $n$,
\begin{thm} \label{thm:RealFiniteDim}
Suppose that $S(\lambda, t_2)$ has one pole $z$ of maximal order
$n$. Then there exists a chain of functions $\{ c_{0*}, c_{1*},
\ldots ,c_{n*} \}$ and $\{ b_0, b_1, \ldots ,b_n \}$ such that
\[ S(\lambda, t_2) = D + \bbmatrix{c_{0*}& c_{1*}& \ldots& c_{(n-1)*}}(\lambda I - \operatorname{Jordan}(z))^{-1}
\bbmatrix{b_0^*\\ b_1^* \\\vdots \\ b_{n-1}^*} \sigma_1
\]
\begin{enumerate}
\item where $c_{0*}$ is a solution of the output differential equation (\ref{eq:OutCC}) with
    the spectral parameter $z$, and $c_{i*}$ is a solution of the differential equation
    \[ z \sigma_{2*} y - \sigma_{1*} \frac{\partial}{\partial t_2}y +
\gamma_* y = \sigma_{2*} c_{i-1*}.
    \]
\item and $b_0$ is a solution of the adjoint output differential equation (\ref{eq:OutCC*}) with the
  spectral parameter $-z^*$, and $b_i$ is a solution of the differential equation
    \[ -z^* \sigma_2 y - \sigma_1 \frac{\partial}{\partial t_2}y +
\gamma y = \sigma_2 b_{i-1}.
    \]
\end{enumerate}
\end{thm}
\textbf{Proof:} Remember that $S(\lambda, t_2)$ is of the form (\ref{eq:SInttw})
\[
S(\lambda, t_2) = \Phi_*(\lambda, t_2^0, t_2) S(\lambda, t_2^0) \Phi^{-1}(\lambda, t_2^0,t_2),
\]
where $\Phi(\lambda, t_2^0, t_2), \Phi_*(\lambda, t_2^0, t_2)$ are
fundamental matrices for the input and the output ODEs,
respectively, and $S(\lambda, t_2^0)$ is a rational in $\lambda$ matrix, for which we can apply the realization theorem:
\[ S(\lambda, t_2^0) = D + C (\lambda I - Jordan[z])^{-1} B
\]
Then we have
\[ \begin{array}{llllllllll}
S(\lambda,t_2^0,t_2 = \Phi_*(\lambda, t_2^0, t_2)  S(\lambda_0, t_2^0) \Phi^{-1}(\lambda, t_2^0,t_2) \\
= \Phi_*(\lambda, t_2^0, t_2) [D + C
\bbmatrix{
\frac{1}{\lambda - z} & \frac{1}{(\lambda - z)^2} & \frac{1}{(\lambda - z)^3} & \ldots & \frac{1}{(\lambda - z)^n} \\
            0         & \frac{1}{\lambda - z} & \frac{1}{(\lambda - z)^2} & \ldots & \frac{1}{(\lambda - z)^{n-1}} \\
            0         &    0                  & \frac{1}{\lambda - z}   & \ldots & \frac{1}{(\lambda - z)^{n-2}} \\
            \vdots    &     \vdots            &    \vdots               & \ddots & \vdots \\
            0         &    0                  &    0                    &  \ldots & \frac{1}{\lambda - z} }
            B] \Phi^{-1}(\lambda, t_2^0,t_2) .
\end{array}
\]
Writing now the fundamental matrices in the Taylour series we obtain that
\[ \begin{array}{llll}
\Phi_*(\lambda,t_2^0, t_2) = \Phi_*(z,t_2^0,t_2) + \sum\limits_{k=1}^\infty \Phi_{k*}(t_2^0, t_2) (\lambda - z)^k =
\sum\limits_{k=0}^\infty \Phi_{i*}(\lambda - z)^k , \\
\Phi^{-1}(\lambda,t_2^0, t_2) = \Phi^{-1}(z, t_2^0, t_2) + \sum\limits_{k'=1}^\infty \Phi_{k'}(t_2^0, t_2)
(\lambda - z)^{k'} = \sum\limits_{k'=0}^\infty \Phi_{k'} (\lambda - z)^{k'}
\end{array} \]
Inserting these expressions for the fundamental matrices, we are able to calculate the Lourent
coefficients explicitly:
\begin{enumerate}
\item The coefficient of $\frac{1}{(\lambda-z)^{n}}$ is
\[ \begin{array}{llllll}
\Phi_{0*} C e_1 e_n^* B \Phi_0 = c_{0*} b_0^* \sigma_1,
\end{array} \]
where we denoted $c_{0*} = \Phi_{0*} C e_1, b_0 = \sigma_1^{-1} \Phi_0 B e_n$.
\item The coefficient of $\frac{1}{(\lambda-z)^{n-1}}$ is
\[ \begin{array}{llllll}
\Phi_{0*} C e_1 e_{n-1}^* B \Phi_0 + \Phi_{0*} C e_2 e_n^* B \Phi_0 +
\Phi_{0*} C e_1 e_n^* B \Phi_1 + \Phi_{1*} C e_1 e_n B \Phi_0 = \\
= c_{0*} e_{n-1}^* B \Phi_0 + \Phi_{0*} C e_2 b_0^*\sigma_1 + c_{0*} e_n^* B \Phi_1 +
\Phi_{1*} C e_1 b_0^* \sigma_1 = \\
= c_{0*}[e_{n-1}^* B \Phi_0 + e_n^* B \Phi_1] +
[\Phi_{0*} C e_2 + \Phi_{1*} C e_1] b_0^* \sigma_1 = \\
= c_{0*} b_1^*\sigma_1 + c_{1*} b_0^* \sigma_1
\end{array} \]
Notice that direct calculations show that $c_{1*}, b_1$ are companion solutions of $c_{0*}, b_0$,
respectively.
\item The coefficient of $\frac{1}{(\lambda-z)^{n-2}}$ is
\[ \begin{array}{llllll}
\Phi_{0*} C e_1 (e_{n-2}^* B \Phi_0 + e_{n-1}^* B \Phi_1 + e_n B \Phi_2) +
 (\Phi_{0*} C e_2 + \Phi_{1*} C e_1)(e_{n-1}^* B \Phi_0 + e_n^* B \Phi_1) + \\
 (\Phi_{0*} C e_3 + \Phi_{1*} C e_2 + \Phi_{2*} C e_1) e_n B \Phi_0 = \\
 = c_{0*} b_2^* \sigma_1 + c_{1*} b_1^* \sigma_1 + c_{2*} b_0^* \sigma_1,
\end{array} \]
where we denote $b_2^* \sigma_1 = e_{n-2}^* B \Phi_0 + e_{n-1}^* B \Phi_1 + e_n B \Phi_2$ and
$c_{2*} = \Phi_{0*} C e_3 + \Phi_{1*} C e_2 + \Phi_{2*} C e_1$. Again, as in the previous case, it is a matter
of simple calculations to show that $c_{2*}, b_2$ are companion solutions to $c_{1*}, b_1$, respectively.
\item For all the other coefficients, by induction, we obtain the desired result.
\end{enumerate}
Notice that $\Phi_{k*}(t_2^0,t_2) = \frac{\partial}{\partial\lambda^k}\Phi_*(\lambda, t_2^0,t_2)|_{\lambda=z},
\Phi_k(t_2^0,t_2) = \frac{\partial}{\partial\lambda^k}\Phi^{-1}(\lambda, t_2^0,t_2)|_{\lambda=z}$ and this is used to
show the connection between $c_{k*}$ and $c_{(k+1)*}$ ($b_k$ and $b_{k+1}$).  \qed

It means that one obtains a vessel with the transfer function
$S(\lambda, t_2)$ by means of the following definitions:
\[ A_1 = \operatorname{Jordan}(z_i), ~~~~A_2(t_2) = 0, ~~~~
C(t_2) = \operatorname{Col}\{ c_{i*}(t_2)\}, ~~~~ B(t_2) =
\operatorname{Row}\{ b _i(t_2)\}
\]
It is a matter of simple calculations to show that all the vessel
conditions are satisfied.

As for the general case, suppose that $A_1 = U \Omega U^{-1}$, where $\Omega$ is the Jordan block form of $A_1$. Then from realization theorem
\ref{tm:NonConsRealiz}, we obtain that
\[ S(\lambda,t_2) = D(t_2) + C(t_2) (\lambda I - A_1)^{-1} B(t_2) \sigma_1 =
D(t_2) + C(t_2) V^{-1} (\lambda I - \Omega)^{-1} V B(t_2) \sigma_1
\]
and one can consider the case where $A_1$ is of the Jordan block form. Suppose that we are given chains of the appropriate sizes
of companion solutions $c^i_{0*}, \ldots, c^i_{(n_i-1)*}$ and $b^i_{0*}, \ldots, b^i_{(n_i-1)*}$ for the output and the adjoint input
ODE with the spectral parameter $\lambda_i$. Suppose that it is given for each eigenvalue $\lambda_i$ of $\Omega$. Then the final vessel is
obtained by defining $C(t_2)$ as a block matrix, where the matrices $c^i_{0*} \ldots c^i_{(n_1-1)*}$ are on the diagonal. The
same construction works for $B(t_2)$, so that the final result and vessel conditions are easily verified.

The construction of the system parameters from the residues of the
given rational matrix function (which one assumes has only simple poles), is known as the Gilbert realization
(see \cite{bib:Kailath}, page 349).

It is also appropriate to emphasize here that we have built in theorem \ref{thm:RealFiniteDim} right and left pole pairs for the matrix function $S(\lambda,t_2)$.
Let us recall first the definitions (from \cite{bib:Inter}). From the formula (\ref{eq:SInttw})
\[ S(\lambda, t_2) = \Phi_*(\lambda,t_2,\tau_2) S(\lambda, t_2^0) \Phi^{-1}(\lambda,t_2,t_2^0)
\]
it follows that the poles of the matrix $S(\lambda, t_2)$ are independent of $t_2$. So, suppose that $z$ is a pole of $S(\lambda,t_2)$, which has the
following Lourent expansion around $z$:
\[ S(\lambda,t_2) = \sum_{j=-q}^\infty (\lambda-z)^j S_j(t_2).
\]
Then $y_0(t_2),\ldots,y_{r-1}(t_2)$ is called a \textit{right pole chain} for $S(\lambda,t_2)$ at $z$, if there exist
additional vectors $y_r(t_2),\ldots,y_{r+q-1}(t_2)$ ($q$ being the order of $z$ as a pole of $S^{-1}(\lambda,t_2)$) such that
$S^{-1}(\lambda,t_2) y(\lambda,t_2)$ is analytic at $z$ with zero of order $r$ at $z$, where
\[ y(\lambda,t_2) = \sum_{j=0}^{r_q-1} (\lambda-z)^j y_j.
\]
The natural number $r$ corresponds to a size of Jordan block (with $z$ on the diagonal) to be constructed shortly.

A \textit{canonical right pole pair} $(A_1,B(t_2))$ is a collection of chains
\[ B(t_2) = y_0^{(1)}(t_2),\ldots,y^{(1)}_{r-1}(t_2);  y_0^{(2)}(t_2),\ldots,y^{(2)}_{r-1}(t_2); \ldots; y_0^{(p)}(t_2),\ldots,y^{(p)}_{r-1}(t_2);
\]
and corresponding Jordan blocks
\[ A_1 = J_1 \oplus J_2 \oplus \cdots \oplus J_p.
\]
with the property above. Moreover, any pair similar to $(A_1,B(t_2))$, i.e., for an invertible matrix $M(t_2)$,
of the form $(M(t_2) A_1 M^{-1}(t_2), M(t_2) B(t_2))$ is called a \textit{right pole pair}.
Analogously, one defines a \textit{canonical left pole pair} $(C(t_2),A_1)$, when $y(\lambda,t_2) S^{-1}(\lambda,t_2)$ is demanded to be
analytic at $z$ with zero of order $r$ at $z$.

In theorem \ref{thm:RealFiniteDim} we have constructed such pairs for each $t_2$, but even more, the left pole pair $(C(t_2),A_1)$ has
an additional property, that $C(t_2)$ satisfies the differential equation in lemma \ref{lm:DEForCB} with the spectral matrix parameter $A_1$.
This in turn means that the left pole chains consist of companion solutions of the output ODE (\ref{eq:OutCC}) for poles of $S(\lambda,t_2)$
as described in theorem \ref{thm:RealFiniteDim}. A similar conclusion is true for $(A_1,B(t_2))$.

Finally, a triple $(C(z,t_2), A_1(z), B(z,t_2)\sigma_1(t_2))$ is called a \textit{pole triple at $z$} of a rational matrix function $S(\lambda,t_2)$ if
$(C(z,t_2),A_1(z))$, $(A_1(z),B(z,t_2)\sigma_1(t_2))$ are the left and the right pole pairs at $z$, respectively. In this case
\[ S(\lambda,t_2) - C(z,t_2) (zI - A_1(z))^{-1} B(z,t_2) \sigma_1(t_2)
\]
is an analytic function at $z$. A (full) \textit{pole triple} of $S(\lambda,t_2)$ is just a direct sum of all local pole triples,
which respects the order of their appearance:
\begin{eqnarray*}
C(t_2) = C(z_1,t_2) \oplus C(z_2,t_2) \oplus \cdots \oplus C(z_p,t_2), \\
A_1 = \operatorname{diag} \{ A_1(z_1), A_1(z_2), \ldots, A_1(z_p) \}, \\
B(t_2)\sigma_1(t_2) = \bbmatrix{B(z_1,t_2) \\ B(z_2,t_2)\\ \vdots \\B(z_p,t_2)} \sigma_1(t_2).
\end{eqnarray*}
A natural question arises, when one can reconstruct $S(\lambda,t_2)\in\mathbf{\mathrm I}$, when one knows its
pole data, i.e., if one knows all pole triples for all $z$'s.

\subsection{Realization theorem (of Mittag-Leffler type)}
Let $S(\lambda, t_2) \in \boldsymbol{\mathcal I}$.
Suppose that in order to solve a Mittag Leffler type problem, we are given a pole triple for each $t_2$, $(X(t_2),T,Y(t_2))$, where we
have denoted the matrix $T$ independent of $t_2$. As we have seen, the pole pairs can be chosen so that the differential equations of lemma
\ref{lm:DEForCB} with spectral matrix parameter $T$ are satisfied. Then the following theorem answers the question of reconstructing
$S(\lambda,t_2)$ from its pole triple.
\begin{thm} \label{thm:ML}
Suppose that we are given a pole triple $(C(t_2),A,B(t_2)\sigma_1)$ for each $t_2$ with constant $A$ and $C(t_2), B(t_2)\sigma_1$ solutions of
the output (\ref{eq:OutCC}) and the inverse input (\ref{eq:InCC*}) ODEs with the matrix spectral parameter $A$.
Suppose also that an analytic function $D(t_2)$ (the value at infinity) is given. Then the matrix function
\begin{equation} \label{eq:SFromPoles}
S(\lambda, t_2) = D(t_2) + C(t_2) (\lambda I - A)^{-1} B(t_2) \sigma_1
\end{equation}
maps solutions of (\ref{eq:InCC}) with spectral parameter $\lambda$ to
solutions of (\ref{eq:OutCC}) with the same spectral parameter iff the linkage conditions (\ref{eq:LinkCond}) are satisfied.
\end{thm}
\textbf{Proof:} The proof is essentially the same as for theorem \ref{tm:NonConsRealiz}. In order to obtain that the function
$S(\lambda,t_2)$ built in (\ref{eq:SFromPoles}) maps solutions of (\ref{eq:InCC}) to solutions of (\ref{eq:OutCC}) with
the same spectral parameter, it is necessary to prove that $S(\lambda, t_2)$ satisfies the following differential equation:
\begin{equation} \label{eq:DS}
 \frac{\partial}{\partial t_2} S(\lambda, t_2) = \sigma_{1*}^{-1} (\sigma_{2*} \lambda + \gamma_*) S(\lambda,t_2) -
S(\lambda,t_2) \sigma_1^{-1} (\sigma_2 \lambda + \gamma),
\end{equation}
i.e., is of the form defined by (\ref{eq:S}).  Differentiating (\ref{eq:SFromPoles}) we obtain
\[\begin{array}{rr}
 \frac{\partial}{\partial t_2} S(\lambda, t_2) & = \frac{\partial}{\partial t_2} D(t_2) +
\sigma_{1*}^{-1} (\sigma_{2*} A + \gamma_*) C(t_2) (\lambda I - A)^{-1}B(t_2)\sigma_1 - \\
& - C(t_2) (\lambda I - A)^{-1} B(t_2) (\sigma_2 A + \gamma).
\end{array} \]
This can be rewritten as
\[ \begin{array}{rrrr}
 \frac{\partial}{\partial t_2} S(\lambda, t_2)  = \frac{\partial}{\partial t_2} D(t_2) -
 \sigma_{1*}^{-1}(t_2) \sigma_{2*} C(t_2) B(t_2) \sigma_1 + C(t_2) B(t_2) \sigma_2(t_2) + \sigma_{1*}^{-1}(t_2) \gamma_*(t_2)D(t_2) + \\
 + \lambda (\sigma_{1*}^{-1}(t_2) \sigma_{2*} D(t_2) - D \sigma_1^{-1}\sigma_2) +
\sigma_{1*}^{-1} (\sigma_{2*} \lambda + \gamma_*)(D(t_2)+ C(t_2) (\lambda I - A)^{-1}B(t_2)\sigma_1(t_2)) - \\
 - (D(t_2) + C(t_2) (\lambda I - A)^{-1} B(t_2)) (\sigma_2 \lambda + \gamma) = \\
=  \frac{\partial}{\partial t_2} D(t_2) -
 \sigma_{1*}^{-1}(t_2) \sigma_{2*} C(t_2) B(t_2) \sigma_1 + C(t_2) B(t_2) \sigma_2(t_2) + \sigma_{1*}^{-1}(t_2) \gamma_*(t_2)D(t_2) + \\
 + \lambda (\sigma_{1*}^{-1}(t_2) \sigma_{2*} D(t_2) - D \sigma_1^{-1}\sigma_2)
 + \lambda (\sigma_{1*}^{-1}(t_2) \sigma_{2*} S(\lambda,t_2) -
 S(\lambda,t_2) \sigma_1^{-1}(\sigma_2 \lambda + \gamma)  .
\end{array}
\]
Demanding further that the differential equation (\ref{eq:DS}) holds for $S(\lambda,t_2)$, we obtain, considering first big $\lambda$ and
then arbitrary one
\begin{enumerate}
    \item $-\sigma_{1*}^{-1} \sigma_{2*} D(t_2) + D(t_2) \sigma_1^{-1} \sigma_2 = 0$ and
    \item $\frac{\partial}{\partial t_2} D(t_2) -\sigma_{1*}^{-1} \gamma_* D(t_2) + D(t_2) \sigma_1^{-1} \gamma +
        \sigma_{1*}^{-1} \sigma_{2*} C(t_2) B(t_2) \sigma_1  - C(t_2) B(t_2) \sigma_2 = 0$.
\end{enumerate}
Denoting $\widetilde D(t_2) = \sigma_{1*} D(t_2) \sigma_1^{-1}$, we obtain that these conditions are
\begin{enumerate}
    \item $- \sigma_{2*} D(t_2) + \widetilde D(t_2) \sigma_2 = 0$ and
    \item $\frac{\partial}{\partial t_2} D(t_2) -\sigma_{1*}^{-1} \gamma_* D(t_2) + \sigma_{1*}^{-1} \widetilde D(t_2) \gamma +
        \sigma_{1*}^{-1} \sigma_{2*} C(t_2) B(t_2) \sigma_1  - C(t_2) B(t_2) \sigma_2 = 0$,
\end{enumerate}
which are exactly the linkage conditions (\ref{eq:LinkCond}). Thus the theorem holds.
\qed


\end{document}